\documentclass[11pt]{article}
\usepackage[pctex32]{graphics}
\usepackage{amssymb}
\usepackage{latexsym}
\usepackage{epsfig}
\usepackage{floatflt}
\usepackage{multirow}
\usepackage{subfigure}
\usepackage{wrapfig}
\usepackage{mathrsfs}
\usepackage[left=1in,top=1in,bottom =1.4in,right=1in,nohead,nofoot]{geometry}
\usepackage{pstcol}

\definecolor{Blue}{rgb}{0.,0.,0.6}
\definecolor{Green}{rgb}{0.,1.,0.0}
\definecolor{Red}{rgb}{1,0.,0.}
\definecolor{Purple}{rgb}{0.5,0.,0.5}

\newcommand{\R}{{\mathbb R}}

\newcommand{\mA}{{\mathsf A}}
\newcommand{\mC}{{\mathsf C}}
\newcommand{\mD}{{\mathsf D}}
\newcommand{\mL}{{\mathsf L}}
\newcommand{\mI}{{\mathsf I}}
\newcommand{\mSigma}{{\mathsf \Sigma}}
\newcommand{\mS}{{\mathsf S}}
\newcommand{\mT}{{\mathsf T}}
\newtheorem{definition}{Definition}[section]
\newtheorem{theorem}[definition]{Theorem}

\begin{document}

\title{Computationally efficient sampling methods for sparsity promoting hierarchical Bayesian models}
\author{D Calvetti \and E Somersalo}
\date{Department of Mathematics, Applied Mathematics and Statistics \\
Case Western Reserve University}

\maketitle
\begin{abstract}
Bayesian hierarchical models have been demonstrated to provide efficient algorithms for finding sparse solutions to ill-posed inverse problems. 
The models comprise typically a conditionally Gaussian prior model for the unknown, augmented by a hyperprior model for the variances. A widely used choice for the hyperprior is a member of the family of generalized gamma distributions. Most of the work in the literature has concentrated on numerical approximation of the maximum a posteriori (MAP) estimates, and less attention has been paid on sampling methods or other means for uncertainty quantification. Sampling from the hierarchical models is challenging mainly for two reasons: The hierarchical models are typically high-dimensional, thus suffering from the curse of dimensionality, and the strong correlation between the unknown of interest and its variance can make sampling rather inefficient. This work addresses mainly the first one of these obstacles. By using a novel reparametrization, it is shown how the posterior distribution can be transformed into one dominated by a Gaussian white noise, allowing sampling by using the preconditioned Crank-Nicholson (pCN) scheme that has been shown to be efficient for sampling from distributions dominated by a Gaussian component. Furthermore, a novel idea for speeding up the pCN in a special case is developed, and the question of how strongly the hierarchical models are concentrated on sparse solutions is addressed in light of a computed example.   
\end{abstract}

\section{Introduction}

A common problem in computational inverse problems is the estimation of an unknown quantity that is a priori believed to be sparse, in the sense that it can be represented with only very few elements of a given basis or frame. In many cases, in particular when the solution is computed numerically, sparsity is replaced by compressibility, meaning that most of the of coefficients in the representation are below a small threshold value. The concepts of sparsity and compressibility are particularly important in the framework of compressed sensing
\cite{Donoho} and sparse dictionary learning \cite{Aharon}. From the very definition of sparsity and compressibility, it is clear that these characterizations are qualitative in nature, as ``most of the coefficients'' is to some extend arbitrary, depends on the dimensionality of the problem, and is open to subjective interpretations. Standard methods for finding sparse solutions include the introduction of sparsity promoting penalties, the most popular being the $\ell^1$-penalty. Thus, if the data $b\in\R^m$ are related to the vector $x\in\R^n$ whose entries are the coefficients of the representation of the unknown in a given basis or frame by $b=f(x)+{\rm noise}$, the standard $\ell^1$-penalized least squares solution is a minimizer of the functional
\[
F_\alpha(x) =  \|b - f(x)\|^2 + \alpha \|x\|_1,
\]
where $\|\,\cdot\,\|$ denotes the Euclidian norm in $\R^m$, $\|\,\cdot\,\|_1$ the $\ell^1$-norm in $\R^n$, and $\alpha>0$ is a regularization parameter. The problem is referred to as basis pursuit \cite{Chen, Candes}, or LASSO \cite{Lasso}, depending on the context. The existence and uniqueness of such solutions depend on the properties of the forward map $f:\R^n\to\R^m$. An alternative, but closely related, approach to sparsity is rooted in the Bayesian analysis of inverse problems \cite{CSnewbook}. In the Bayesian framework, the prior belief of the sparsity of the signal is encoded into a prior that favors sparse solutions, and single point estimates such as the maximum a posteriori (MAP) or posterior mean estimates are generated in the hope that they have the desired sparsity properties \cite{CSgaussian,L2Magic}. In this article, we restrict our attention to a particular family of Bayesian sparsity promoting priors, namely hierarchical Gaussian priors augmented with a hyperprior from the family of generalized gamma distributions \cite{Analysis}, reviewed in section~\ref{sec:hierarchical}. In that section, we also give a brief review of an algorithm to compute the MAP estimates, the iterative alternating sequential (IAS) algorithm, which has been shown to generate sparse, or compressible solutions at a relatively small computational cost.  Under certain conditions, the IAS algorithm has been demonstrated to converge to a unique solution that approximates the $\ell^1$-penalized solution \cite{L2Magic}. The convergence properties are leveraged in a hybrid algorithm that combines different choices of generalized gamma priors \cite{Hybrid}. For an alternative but related way of estimating the hyperparameters based on data, we refer to \cite{Bortoli,Vidal}.

A pertinent question concerning the IAS solution or other MAP estimates and, in general, algorithms searching for a single estimate, is how representative the sparse solution is. A common criticism of the MAP estimate is that it may capture poorly the posterior distribution and could be unstable with respect to perturbations in the data, in particular in the case when the posterior density is multimodal, whereas the posterior mean estimate, when calculated by means of Monte Carlo sampling from the posterior represents a more reliable alternative as a single point estimator. The question is closely related to the wider question of uncertainty quantification under sparsity constraints, often addressed by  Markov chain Monte Carlo  (MCMC) sampling \cite{Kaipio,Lucka}. While the use of MCMC methods to explore posterior densities is rather standard, hierarchical sparsity promoting models are known to pose significant challenges to sampling methods.  The problems are twofold: The problems are typically high-dimensional, and the typically strong correlation between the unknown of primary interest and the hyperparameters  lead to poor mixing and extremely slow convergence of the samplers.
Remedies to the latter problem have been proposed in the literature, including appropriate changes of variables, see, e.g. \cite{Agapiou,Betancourt,Roberts}.  For contributions to quantify the uncertainties in the inverse problems based on hierarchical prior models, we refer to \cite{Agapiou,Agrawal,Glaubitz}.

In this article, we propose changes of variables specifically tailored for the hierarchical Bayesian models that the IAS algorithm is based on. The main goal of this study is to address the problems arising from the high dimensionality of the problem. The change of variables that we propose, combined with the preconditioned Crank-Nicholson (pCN) sampling algorithm \cite{Cotter}, leads to an easy to implement sampling algorithm that is fast and relatively easy to tune. The algorithm provides  and efficient sampling of the posterior, in particular when the hypermodel is based on a gamma distribution or generalized gamma distributions near it. Hypermodels with strong non-convexity remain a challenge, although for the inverse gamma hyperprior, we propose a particular radial  parametrization generalizing the pCN sampler  that improves significantly the convergence. Another novel contribution of this work is an automatic way to choose the hypermodel parameters in the hybrid IAS algorithm, thus completing the work in \cite{Hybrid}. The computed examples using this sampler to explore posterior with generalized gamma hyperpriors seem to suggest a rather surprising result, namely that while the MAP estimate is itself very definitely compressible, neither the Monte Carlo samples around it nor the corresponding posterior mean are necessarily sparse or compressible unless the hypermodel is chosen to be strongly non-convex. Therefore, one can argue that the MAP estimates based on hierarchical Gaussian models are optimal when capturing the sparse nature of the unknown is important.

\section{Hierarchical models and sparsity}\label{sec:hierarchical}

Consider the inverse problem of estimating an unknown $x\in\R^n$ from noisy observations of a linear transformation of it,
\[
 b = f(x) + e,
\]
where $f:\R^n\to\R^m $ is assumed to be a known function. In the Bayesian framework, all parameters not known exactly are modeled as random variables. In the rest of the paper, we denote random variables by capital letters, and their realizations by lowercase letters. The stochastic extension of the observation model is therefore
\begin{equation}\label{obs model}
 B = f\big(X\big) + E,
\end{equation}
where it is assumed that $X$ and $E$ are mutually independent random variables. If the probability distribution of the noise $E$ is given in terms of a probability density function $\pi_E$, the likelihood model for $B$ is 
\[
 \pi_{B\mid X}(b\mid x) = \pi_{E}(b - f(x)).
\]  
If we assume that $E$ is a zero mean Gaussian noise with positive definite covariance matrix $\mSigma\in\R^{m\times m}$, we obtain the
likelihood model 
 \[
 \pi_{B\mid X}(b\mid x) \propto {\rm exp}\left(-\frac 12 (b - f(x))^\mT \mSigma^{-1}(b - f(x))\right).
\]
Furthermore, by a standard whitening argument, multiplying both $b$ and $f$ by a symmetric factor $\mS$ of the precision matrix, $\mSigma^{-1} = \mS^\mT \mS$, we may assume without loss of generality that $\mSigma = \mI_m$, the $m\times m$ identity matrix, and the likelihood model simplifies to
\begin{equation}\label{lkh}
 \pi_{B\mid X}(b\mid x) \propto {\rm exp}\left(-\frac 12 \|b - f(x)\|^2 \right).
\end{equation}

To express the sparsity or compressibility belief, we consider a conditionally Gaussian prior model,
\[
 \pi_{X\mid\Theta}(x\mid\theta) =\left(\frac 1{(2\pi)^n\theta_1 \cdots \theta_n}\right)^{1/2}{\rm exp}\left( -\frac 12 \sum_{j=1}^n \frac{x_j^2}{\theta_j}\right) \propto {\rm exp}\left( -\frac 12 \sum_{j=1}^n \frac{x_j^2}{\theta_j} -\frac 12 \sum_{j=1}^n \log\theta_j \right).
\] 
Furthermore, we assume that the variances $\Theta_j$ are mutually independent and distributed according to the generalized gamma distribution,
\[
 \Theta_j \sim \pi_{\Theta_j}(\theta_j \mid \vartheta_j,\beta,r)  = \frac{|r|}{\Gamma(\beta)\vartheta_j}\left(\frac{\theta_j}{\vartheta_j}\right)^{r\beta -1}
 {\rm exp}\left( - \left(\frac{\theta_j}{\vartheta_j}\right)^r\right), \quad 1\leq j\leq n,
\]
where $r\neq 0$ and the shape parameter $\beta>0$ is assumed to be the same for all $j$ while the scale parameter $\vartheta_j$ may differ for every $j$. Taking into account the mutual independency of the variances, we can write the joint prior model in the form
\begin{eqnarray}\label{joint prior}
 \pi_{X,\Theta}(x,\theta) &=& \pi_{X\mid\Theta}(x\mid\theta)\pi_{\Theta}(\theta\mid \vartheta,\beta,r) \nonumber \\
 &\propto& {\rm exp}\left( -\frac 12 \sum_{j=1}^n \frac{x_j^2}{\theta_j}  - \sum_{j=1}^n\left(\frac{\theta_j}{\vartheta_j}\right)^r +\left(r\beta - \frac 32\right) \sum_{j=1}^n \log\frac{\theta_j}{\vartheta_j} \right).
 \end{eqnarray}

Next we proceed by {\em non-dimensionalizing} the model. Introduce the non-dimensional variables
\begin{equation}\label{scaling}
 \xi_j = \frac{x_j}{\sqrt{\vartheta_j}}, \quad \lambda_j = \frac{\theta_j}{\vartheta_j},
\end{equation}
and observe that after the change of variables, the components $(\xi_j,\lambda_j)$ are a priori independent and identically distributed (i.i.d.),
\[
 (\xi_j,\lambda_j) \sim  {\rm exp}\left( -\frac 12 \frac{\xi_j^2}{\lambda_j}  - \lambda_j^r  +\left(r\beta - \frac 32\right)  \log\lambda_j \right),
\]
and the posterior density can be written as
 \begin{equation}\label{scaled posterior}
 \pi_{\Xi,\Lambda\mid B}(\xi,\lambda \mid b) \propto
  {\rm exp}\left( -\frac 12 \|b - f\big(\mD_{\vartheta}^{1/2} \xi\big) \|^2 -\frac 12 \sum_{j=1}^n \frac{\xi_j^2}{\lambda_j}  - \sum_{j=1}^n\lambda_j^r  +\left(r\beta - \frac 32\right) \sum_{j=1}^n \log\lambda_j \right),
\end{equation}
where $\mD_\vartheta\in\R^{n\times n}$ is a diagonal matrix with the vector $\vartheta$ along its diagonal.  Observe that in the case of a linear forward model, $f(x) = \mA x$, where $\mA\in\R^{m\times n}$, we have 
\begin{equation}\label{linear}
f(\mD_{\vartheta}^{1/2}\xi) = \mA \mD_{\vartheta}^{1/2} \xi,
\end{equation}
where the transformation $\mA\mapsto \mA \mD_{\vartheta}^{1/2}$ is a column scaling of the forward map. In \cite{BTOP,L2Magic}, it was demonstrated that this scaling, which is part of the prior, can be associated to the sensitivity of the unknowns $x_j$ to the data. For completeness, the relevant results are stated here without proof.

Given an observation model (\ref{obs model}), the signal-to-noise ratio (SNR) is defined as
\[
 {\rm SNR} = \frac{{\rm E}(\|B\|^2)}{{\rm E}(\|E\|^2)} =\frac{{\mbox{signal power}}}{{\mbox{noise power}}}.
\] 
Furthermore, assume that the variable $X$ is supported in a set $S\subset\{1,2,\ldots,n\}$. When restricting the support of $X$ to $S$, we denote the corresponding SNR by $SNR_S$. 

\begin{definition}
The random variable satisfies the SNR-exchangeability condition with respect to the observation model, if
\[
 SNR_S = SNR_{S'} \mbox{ for all subset $S,S'\in\{1,2,\ldots,n\}$ with the same cardinality.}
 \]
\end{definition}

Recall that in probability theory, exchangeability is defined by the condition that the probability distribution of a multivariate random variable is invariant under permutation of the variables. Exchangeability of the prior density (before non-dimensionalization) seems a very natural requirement in the Bayesian setting, as it states that all variables are treated with no  a priori preferences, however, it is known that exchangeable priors often lead to useless reconstructions, e.g., by favoring sources close to receivers in geophysical or medical imaging applications. The SNR-exchangebility provides a natural resolution of this dilemma. 
Exchangeability implies the weaker SNR-exchangebility, while converse is not true. The following theorem, however, shows that while the SNR-exchangeability does not imply exchangeability, it does, in the linear case, guarantee that each component of $x$ has an equal chance to explain the data.

\begin{theorem}
Assume that the random variable $X$ satisfies the SNR-exchangeability condition with respect to the linear observation model (\ref{obs model}) with (\ref{linear}), and the prior is the hierarchical prior (\ref{joint prior}). Then the scaling parameters must satisfy
\[
 \vartheta_j = \frac{C}{\|a^{(j)}\|^2},
\]
where $a^{(j)}\in\R^n$ is the $j$th column of the matrix $\mA$, and the constant $C$ depends on the parameters $r$ and $\beta$ as well as the SNR.
\end{theorem}

In \cite{Analysis}, an explicit formula for $C$ is given.  Observe  that with the assumption of the SNR-exchangeability, the column scaling by $\mD_\vartheta^{1/2}$ corresponds to the scaling of the columns of $\mA$ by the column norms. This is a common procedure in optimization \cite{DennisSchnabel} for balancing, and in geophysics and medical imaging, where the scaling is related to sensitivity analysis, as the sensitivity of the data to the component $x_j$ in the linear case is given by
\[
 s_j = \|\partial_{x_j} \mA x\| = \|a^{(j)}\|,
\] 
hence effectively the model is scaled to balance the sensitivity. This point will becomes more clear when we discuss the MAP estimates and hypermodels.

\section{Exploring the posterior density}

In this section, after a brief review the IAS algorithm for estimating the maximum a posteriori (MAP) estimate, we proceed to discuss the exploration of the posterior density using Markov chain Monte Carlo (MCMC) methods.

\subsection{MAP estimate and the IAS algorithm}

The MAP estimate corresponding to the posterior density (\ref{scaled posterior}) is, by definition, 
\[
 (\xi_{\rm MAP},\lambda_{\rm MAP}) = {\rm argmax}\{\pi_{\Xi,\Lambda}(\xi,\lambda\mid b)\},
\]
provided that such maximizer exists. Equivalently, the MAP estimator, if it exists, is also a minimizer of the Gibbs energy,
\begin{equation}\label{gibbs}
{\mathscr E}(\xi,\lambda) = \frac 12 \|b - f\big(\mD_{\vartheta}^{1/2} \xi\big) \|^2 +\frac 12 \sum_{j=1}^n \frac{\xi_j^2}{\lambda_j}  + \sum_{j=1}^n\lambda_j^r  - \left(r\beta - \frac 32\right) \sum_{j=1}^n \log\lambda_j.
\end{equation}
 The IAS algorithm searches for the MAP estimate through the alternating steps as follows. Given an initial $\lambda^0$, set the counter at $t=1$, and iterate the steps until a convergence criterion is met,
 \begin{itemize}
 \item[(a)] Update $\xi$ by defining
 \[
  \xi^t = {\rm argmin} \{{\mathscr E}(\xi,\lambda^{t-1})\} = {\rm argmin}\left\{\frac 12 \|b - f\big(\mD_{\vartheta}^{1/2} \xi\big) \|^2 +\frac 12 \sum_{j=1}^n \frac{\xi_j^2}{\lambda_j^{t-1}} \right\},
\]
 \item[(b)] Update $\lambda$ by defining
\[
  \lambda^t = {\rm argmin} \{{\mathscr E}(\xi^t,\lambda)\} = {\rm argmin}\left\{\frac 12 \sum_{j=1}^n \frac{(\xi_j^t)^2}{\lambda_j}  + \sum_{j=1}^n\lambda_j^r  - \left(r\beta - \frac 32\right) \sum_{j=1}^n \log\lambda_j\right\},
\] 
\item[(c)] Advance the counter, $t\rightarrow t+1$ and check for convergence. 
 \end{itemize} 
 Observe that if the forward model is linear, the first step is a standard least squares problem, and regardless of the linear model, the second step is a component-wise updating problem requiring the solution of the first order optimality condition
 \begin{equation}\label{optimality}
 \frac{\partial}{\partial\lambda_j}\left(\frac{(\xi_j^t)^2}{\lambda_j}  + \lambda_j^r  - \left(r\beta - \frac 32\right)  \log\lambda_j\right) = 0.
\end{equation}
For certain values of $r$, (\ref{optimality}) admits a closed form solution. When $r=1$ and $\eta = \beta- 3/2>0$, corresponding to gamma hyperpriors, the updating formula for $\lambda^t$ becomes
 \[
  \lambda_j^t = \frac 12\left(\eta + \sqrt{\eta^2 + 2 (\xi_j^t)^2}\right), \quad \eta = \beta - \frac 32,
\]
while, when $r=-1$, corresponding to inverse gamma hyperpriors, we have
\begin{equation}\label{lambda inv gamma}
 \lambda_j^t = \frac 1\kappa\left(\frac{(\xi_j^t)^2}2 + 1\right), \quad \kappa = \beta + \frac 32.
\end{equation}
For choices of $r$ that do not admit closed form solutions, the values $\lambda^t_j$ satisfy
\[
 \lambda_j^t = \varphi\big(|\xi_j^t|\big),
\]
where the function $\varphi:\R_+\to\R$ solves the non-linear initial value problem
\begin{equation}\label{varphi}
 \varphi'(t) = \frac{2 t\varphi(t)}{2r^2\varphi(t)^{r+1} + t^2},\quad \varphi(0) = \left(\beta - \frac{3}{2r}\right)^{1/r},
\end{equation}
obtained from implicit differentiation of the equation (\ref{optimality}).
The evaluation of the function $\varphi$ at the points $|\xi_j^t|$ can be done efficiently by sorting the values $|\xi_j^t|$ in ascending order and solving the initial value problem by any standard numerical integrator at those points. For details, we refer to \cite{Analysis}.

The existence and uniqueness of the MAP estimate, as well as the convergence of the IAS algorithm are not obvious. The following result has been proved in \cite{L2Magic}  for the case when $r=1$, $\beta> 3/2$, and $\mA$ is linear.

\begin{theorem}\label{th:eta}
Assume that the forward map is linear, $r=1$ and $\eta = \beta - 3/2>0$. Then the Gibbs energy functional (\ref{gibbs}) has a unique minimizer $(\xi^*,\lambda^*)$, and the IAS algorithm converges to that minimizer. Moreover, it holds that 
\[
 \lambda_j = \frac 12\left(\eta + \sqrt{\eta^2 + 2 (\xi^*_j)^2}\right) = f_j(\xi^*_j),
\]
and $\xi_j^*$ is the unique minimizer of the functional
\[
 \widehat{\mathscr E}(\xi) = {\mathscr E}(\xi, f(\xi)).
\]
Moreover, 
\[
 \lim_{\eta\to 0+} \widehat {\mathcal E}(\xi) = \frac 12 \|b -  \mA\mD_\vartheta^{1/2} \xi\|^2 + \sqrt{2} \sum_{j=1}^n |\xi_j|,
\]
and the IAS solution converges to the minimizer of the above right hand side.   
\end{theorem}
 
The above theorem underlines the role of the parameter $\eta$ in promoting sparsity of the solution. For other hyperparameter values, in particular for $r<1$, the uniqueness of the minimizer may not be guaranteed even in the linear case, and in fact, in some cases the algorithm is known to converge to local minima. The sparsity promoting nature of the prior can often be understood by restricting the Gibbs energy functional to the manifold defined by the minimization condition (b) for $\lambda$. For instance, in the case of the inverse gamma model $r=-1$, defining $g_j(\xi_j)$ through the formula (\ref{lambda inv gamma}), we obtain 
\[
 {\mathscr E}(\xi,g(\xi)) = \frac 12 \|b - f\big(\mD_{\vartheta}^{1/2} \xi\big) \|^2  +\sum_{j=1}^n \log\left(1+\frac{\xi_j^2}{2} \right)^\kappa + {\rm constant}, \quad \kappa = \beta +\frac32
\]
corresponding to a prior of the form
\[
\pi_\Xi(\xi) \propto \prod_{j=1}^n \frac{1}{(1+\xi_j^2/2)^\kappa},
\]
which is strongly sparsity promoting. At the limit $\beta \to 0+$, $\kappa \to \frac32$, this prior converges to the Student distribution with $\nu=2$; see \cite{Analysis} for details.

In \cite{Hybrid} a hybrid IAS algorithm was proposed and investigated. The idea behind the hybrid model is to run first the IAS algorithm by using the parameter value $r=1$, that  is guaranteed to converge in the linear case, then switch to a greedier  scheme by choosing a generalized gamma hyperprior with $r<1$.  
When the switching occurs, the two hypermodels are matched using the following two criteria:
\begin{enumerate}
\item Whenever $x_j=0$, we require that the baseline values for $\theta_j$ coincide. This way, the a priori variance of the background outside the support of $x$ is consistently defined independently of the model. 
\item The marginal expected value for $\theta_j$ is equal using both models.
\end{enumerate}
  
Let us denote by $(r_1,\beta_1,\varphi_1)$ and $(r_2,\beta_2,\varphi_2)$ the hyperparameter values for two models, where in practice $r_1=1$. 
From the initial condition in (\ref{varphi}), we conclude that the compatibility condition implies that
\begin{equation}\label{compatibility}
 \vartheta_1\left(\beta_1 - \frac{3}{2r_1}\right)^{1/r_1} =  \vartheta_2\left(\beta_2 - \frac{3}{2r_2}\right)^{1/r_2}.
\end{equation}
Recalling the expectation of generalized gamma distribution, the second condition can be written as
\begin{equation}\label{compatibility2}
\vartheta_1\frac{\Gamma(\beta_1 + \frac 1{r_1})}{\Gamma(\beta_1)} = \vartheta_2\frac{\Gamma(\beta_2 + \frac 1{r_2})}{\Gamma(\beta_2)}.
\end{equation} 
The assumed finiteness of the expectation poses restrictions to possible parameter values.
In the section of computed examples, these conditions are discussed in detail for several special cases. We point out that the second condition was not considered in the cited paper, and it is introduced here for the first time to make the hyperparameter selection automatic.

In the following, we shall use the IAS algorithm, and the hybrid scheme in particular, to find an appropriate initial point for the MCMC sampling.

\subsection{Sampling with a Gaussian prior: Preconditioned Crank-Nicholson}

In preparation of the reparametrization of the hypermodels, we recall some known results concerning random draws from Gaussian distributions. Let $\mC\in\R^{n\times n}$ be a symmetric positive definite matrix. Assuming that a symmetric factorization  $\mC = \mL^\mT\mL,$ such as the Cholesky factorization of the matrix, is available, independent random draws from the normal distribution ${\mathcal N}(x\mid 0,\mC)$ can be generated through the formula
\[
 X = \mL^\mT W, \quad W\sim{\mathcal N}(0,\mI_n),
\] 
where $\mI_n$ is the $n\times n$ identity matrix.
The independent sampling, while generating a sequence $X^j$ of independent draws from the distribution, does not provide any way to control the step size $\|X^j - X^{j-1}\|$. Step size control is fundamental when the Gaussian distribution is used as a proposal distribution for exploring posterior distributions. A way to enable step size control is to consider the sequence
\[
 X^j = \sqrt{1-h^2}X^{j-1} + h \mL^\mT W, \quad W\sim{\mathcal N}(0,\mI_n),\quad X^0\sim{\mathcal N}(0,\mC),
\] 
where $0<h<1$.  It is a straightforward matter to check by induction that $X^j$ is Gaussian zero mean random variable, with covariance equal to $\mC$, therefore, the produced sequence is distributed according to ${\mathcal N}(0,\mC)$. The parameter $h$ controls the step size in the sense that
\begin{eqnarray*}
 {\rm E}\|X^j - X^{j-1}\|^2 &=& {\rm trace}\big({\rm E}(X^j-X^{j-1}) (X^j-X^{j-1})^\mT\big) \\
  &=& 2(1 - \sqrt{1-h^2})\,{\rm trace}(\mC) \approx h^2\,{\rm trace}(\mC)
\end{eqnarray*} 
for $h$ small. The disadvantage, compared to independent sampling, corresponding to $h=1$,  is that consecutive samples are correlated, since
\[
 {\rm E}(X^{j-1}(X^j)^\mT) =  \sqrt{1-h^2}\mC.
\] 
These results are well known in the literature. In  \cite{LiuWest}, this observation was used to modify Gaussian mixtures so as to avoid artificial diffusion in mixture model. In \cite{Cotter}, this observation is used to define a Metropolis-Hastings type algorithm that is efficient for high dimensional inverse problems, known as {\em preconditioned Crank-Nicholson} (pCN) algorithm, which is a key tool in this paper. In the cited article, the authors consider distributions of the type
\[
 \pi(X) \propto e^{-\Phi(x)} {\mathcal N}(x\mid 0,\mC).
\] 
By defining a proposal $y$ drawn from a non-symmetric proposal ${\mathcal N}(y\mid \sqrt{1-h^2}\, x^{j-1}, h^2\mC)$,  it is easy to check that the Metropolis-Hastings acceptance ratio reduces to
\[
 \alpha(x^{j-1},y) = e^{\Phi(x^{j-1}) - \Phi(y)},
\]
that is, the potentially high-dimensional Gaussian distribution does not appear in the acceptance ratio, thus avoiding the problems that in the high dimension limit lead to practically automatic rejection, in line with the Cameron-Martin theorem.

\subsection{Reparametrization}

We are now ready to introduce a particular reparametrization that makes it possible to take advantage of the pCN algorithm for MCMC sampling hierarchical models.  Consider the posterior density (\ref{scaled posterior}) and let $\R^{2n} = \otimes_{j=1}^N\R^2$. In each subspace $\R^2$ introduce the new parameters $(v_j,\tau_j)$ such that
\[
 v_j^2 = \frac{\xi_j^2}{\lambda_j}, \quad \tau_j^2 = 2 \lambda_j^r, \quad \lambda_j>0.
\] 
In the transformation of the old variables $(\xi_j,\lambda_j)$ in terms of the new ones $(v_j,\tau_j)$, we choose the signs so that
\[
\lambda_j = \left(\frac{\tau_j^2}2\right)^{1/r}, \quad \xi_j = v_j \sqrt{\lambda_j} = v_j \left(\frac{\tau_j^2}2\right)^{1/2r} =\frac{v_j  |\tau_j|^{1/r}}{2^{1/2r}}.
\]
The reparameterization in the probability density requires the determinant of the Jacobian in each subspace $\R^2$, 
\[
J(v_j,\tau_j) =  {\rm det}\left(\left[\begin{array}{cc}\displaystyle{\frac{\partial \lambda_j}{\partial \tau_j}} & \displaystyle{\frac{\partial \lambda_j}{\partial v_j} } \\
 \displaystyle{\frac{\partial \xi_j}{\partial \tau_j} }& \displaystyle{\frac{\partial \xi_j}{\partial v_j}}\end{array}\right]\right)
 = {\rm det}\left(\left[\begin{array}{cc}  \displaystyle{\frac{2}{r 2^{1/r}}\frac{|\tau_j|^{2/r}}{\tau_j}} & 0 \\ \displaystyle{\frac{1}{r 2^{1/2r}}\frac{v_j |\tau_j|^{1/r}}{\tau_j} }
 & \displaystyle{\frac{ |\tau_j|^{1/r}}{2^{1/2r}}}\end{array}\right]\right) =\frac{2^{1-3r/2}}{r} \frac{|\tau_j|^{3/r}}{\tau_j},
\] 
The posterior density, written in terms of the new variables $(v,\tau)$ then becomes
\begin{eqnarray*}
&&  \pi_{V,T}(v,\tau\mid b)  =  \prod_{j=1}^n |J(v_j,\tau_j)| \pi_{\Xi,\Lambda}(\xi(v,\tau),\lambda(v,\tau)\mid b) \\
 &&\propto
   {\rm exp}\bigg(-\frac 12 \left\|b -  f\left(\frac1{2^{1/2r}} \mD_\vartheta^{1/2} (v |\tau|^{1/r})\right) \right\|^2 + \left(r\beta-\frac 32\right)\sum_{j=1}^n \log|\tau_j|^{2/r} \\
   &&\phantom{XXXXXXXXXXXXXXXXXXXX}
  +\sum_{j=1}^n \log |\tau_j|^{3/r-1}-\frac 12\|v\|^2 - \frac 12 \|\tau\|^2 \bigg) \\
   && \propto e^{-\Phi(v,\tau)}{\mathcal N}(v,\tau\mid 0,\mI_{2n}),
 \end{eqnarray*} 
 where the products are understood to be component-wise, i.e.,
 \[
  \big(v|\tau|^{1/r}\big)_j =v_j |\tau_j|^{1/r}, \quad 1\leq j\leq n,
\]
and the potential function $\Phi$ is defined as   
 \[
  \Phi(v,\tau) = \frac 12 \left\|b -  f\left(\frac1{2^{1/2r}} \mD_\vartheta^{1/2} (v |\tau|^{1/r})\right) \right\|^2 - \left(2\beta-1\right)\sum_{j=1}^n \log |\tau_j| .
\]  
In the following, we will consider four special cases: When $r=1$, the hyperprior is the gamma distribution, and we have
\[
r=1:\quad   \Phi(v,\tau) = \frac 12 \left\|b -  f\left(\frac1{\sqrt{2}} \mD_\vartheta^{1/2} (v |\tau|)\right) \right\|^2 - \left(2\beta -1\right)\sum_{j=1}^n \log |\tau_j| ,
\]  
while if $r=-1$, the hyperprior is the inverse gamma distribution, and 
\[
 r=-1:\quad\Phi(v,\tau) =  \frac 12 \left\|b -  f\left({\sqrt{2}} \mD_\vartheta^{1/2} \frac v{ |\tau|}\right) \right\|^2 - \left(2\beta-1\right)\sum_{j=1}^n \log |\tau_j| .
\] 
In the case $r=1/2$, we have
\[
 r = \frac 12: \quad   \Phi(v,\tau) = \frac 12 \left\|b -  f\left(\frac1{2} \mD_\vartheta^{1/2} (v  \tau^{2})\right) \right\|^2 - \left(2\beta -1\right)\sum_{j=1}^n \log |\tau_j| ,
\] 
and for $r = -1/2$, we have
\[
 r = -\frac 12: \quad   \Phi(v,\tau) = \frac 12 \left\|b -  f\left( 2 \mD_\vartheta^{1/2} (\frac v {\tau^{2}})\right) \right\|^2 - \left(2\beta -1\right)\sum_{j=1}^n \log |\tau_j| ,
\] 
In the following section, we use the proposed reparametrization to explore the posterior densities corresponding to these four choices of $r$.

\section{Computed examples}

In this section, we investigate numerically the effectiveness of the reparametrization of the problem in combination with the pCN algorithm, with special emphasis on the role of the parameter $r$.  We start by discussing the model problems used in testing the proposed sampler.

\subsection{Model problem}

We  consider here a linear inverse problem in the form of a one-dimensional deconvolution. Let $g:[0,1]\to\R$ be the function to be estimated from noisy observations of the convolution of the signal with a Gaussian kernel,
\begin{equation}\label{forward model}
 b_j = \int_0^1 a(t_j - s) g(s) ds + \varepsilon_j \quad a(t) = A\,{\rm exp}\left(-\frac 1{2 w^2}t^2\right), \quad 0<t_1<\ldots<t_m<1,
\end{equation}
of width $w=0.02$ and amplitude $A=6.2$. We discretize the convolution integral using a piecewise constant approximation with $n = 128$  intervals, and assume that the observation points $t_j$ coincide with every sixth discretization node $s_k$, that is, $t_j  = s_{1+ 6(j-1)}$, $1\leq j\leq m = 22$. This yields the approximation
\[
 b = \mA z + \varepsilon, \quad \mA\in\R^{m\times n}, \quad z_j = g(s_j) .
\] 
To generate the data, we assume that the generative model $g$ is a piecewise constant function. 
To avoid the inverse crime of using the same model for data generation and for the solution of the inverse problem, the data are generated by using a fine discretization mesh with 1000 discretization intervals, and subsequently Gaussian noise, $\varepsilon$ from ${\mathcal N}(0,\sigma^2\mI_m)$ with $\sigma = 0.03$ is added to them. The generative model and the noisy data are shown in the top left panel of  Figure~\ref{fig:IAS results}.

To find a sparse representation of the discrete signal $z$, let $\mL\in\R^{n \times n}$ be the finite difference matrix,
\[
 \mL = \left[\begin{array}{rrrr} 1 & & &\\ -1 & 1 & & \\ & \ddots & \ddots & \\  & & -1&\;1\end{array}\right],
\] 
and express $z$ in terms of $x\in\R^n$ as $\mL z = x$. Implicitly, we assume here a boundary condition $z_0 = g(0) = 0$. Therefore, we can write the forward model as
\[
 b = \mA\mL^{-1} x + \varepsilon,
\]
and after scaling the data and the forward map by $1/\sigma$,
\[
 \widehat b = \frac 1\sigma b, \quad \widehat \mA = \frac 1\sigma \mA \mL^{-1},
\]
thus whitening the noise, we arrive at the expression of the problem in standard form,
\begin{equation}\label{linear}
 \widehat b = \widehat\mA x + e, \quad e\sim{\mathcal N}(0,\mI_m).
\end{equation}

We point out that the vanishing boundary value of the unknown $z$ at the endpoint $t=0$ makes the data non-uniformly sensitive to the components of the vector $x$. The sensitivity of the data to the components of the vector $x$ is not addressed here: we refer to \cite{L2Magic} for the discussion of the topic.

\subsection{Hypermodels and MAP estimates}

We test the sampling algorithm with four different generalized gamma hypermodels, corresponding to  $r_1=1$, $r_2=1/2$, $r_3=-1/2$ and $r_4=-1$.  We denote the corresponding hyperparameter values by $(\beta_j,\vartheta_j)$, $1\leq j\leq 4$.  Here, for simplicity, we choose all components of the vectors $\vartheta_j$ equal, so this parameter can be treated as a scalar.

In order to initialize the MCMC algorithm without the need of a long burn-in run, we first compute a MAP estimate for each case using the hybrid IAS algorithm: 
\begin{enumerate}
\item Phase I: Run the IAS algorithm with values $(r_1,\beta_1,\vartheta_1)$  until convergence criterion is met;
\item Phase II: If $j  > 1$, continue the IAS iteration with values  $(r_j,\beta_j,\vartheta_j)$ until convergence criterion is met.
\end{enumerate}
The compatibility conditions (\ref{compatibility})-(\ref{compatibility2}) give an automatic way to set the hyperparameters for $j>1$. Recalling that by Theorem~\ref{th:eta} for $r=1$,  a value $\beta_1$ close to $3/2$ promotes sparsity, we write $\beta_1 = 3/2 +\eta$, where $\eta>0$ is small.  We also set the value $\vartheta_1$, allowing automatic selection of the hyperparameters for $j>1$. The compatibility conditions with $r_1=1$ yield
\[
 \vartheta_j\left(\beta_j - \frac{3}{2r_j}\right)^{1/r_j} = \vartheta_1 \eta, \quad
  \vartheta_j \frac{\Gamma(\beta_j + \frac 1{r_j})}{\Gamma(\beta_j)} = \vartheta_1\left(\eta+\frac32\right),\quad j=2,3,4.
\]
Straightforward calculations based on the properties of the gamma function lead to the following formulas for the parameters:

{\bf Case} $r_2 = 1/2$: We have
\[
 \beta_2  = \frac{6m+1  + \sqrt{48 m +1}}{2(m-1)}, \quad  \mbox{where } m = 1+\frac3{2\eta},
\]
and
\[
 \vartheta_2 =  \vartheta_1 \frac{\eta}{(\beta_2 -3)^2}.
\] 

{\bf Case} $r_3 = -1/2$: We have
\[
 \beta_3  = \frac{6+3m \pm \sqrt{m^2+80 m}}{2(m-1)}, \quad  \mbox{where } m = 1+\frac3{2\eta},
\]
and
\[
 \vartheta_3 =  \vartheta_1 \eta{(\beta_2 +3)^2}.
\] 

{\bf Case} $r_3 = -1$: We have
\[
 \beta_4  = 1 + \frac{5}{3}\eta,
\]
and
\[
 \vartheta_4 =  \vartheta_1 \eta \left(\beta_4 +\frac 32\right).
\] 
The numerical values used in the computations are give in Table~\ref{tab:params}.

In Phase II, the IAS iterations start with the final $\theta$ of Phase I. In both phases, the IAS iterations stop as soon as
\[
 \frac{\|\theta^{t-1} - \theta^t\|}{\|\theta^{t-1}\|} < 0.005.
\]

\begin{table}
\centerline{
\begin{tabular}{l|llll}
 & $r=1$ & $r=1/2$ & $r=-1/2$ & $r=-1$ \\
 \hline
 $\beta$ & 1.501 & 3.0918 &  2.0165  & 1.0017\\
 $\vartheta$ & $5\times 10^{-2}$ &$5.9323\times 10^{-3}$  & $1.2583\times 10^{-3}$ & $1.2308\times 10^{-4}$
\end{tabular}
}
\caption{\label{tab:params} Hyperparameter values used in the computations. The values $\vartheta$ for the hybrid models are determined by requiring that at $x_j=0$, the values $\theta_j$ given by the IAS algorithm are independent of the hypermodel, and that the marginal expectations for $\theta_j$ are independent of the model.} 
\end{table}

\begin{figure}[ht!]
\centerline{\includegraphics[width=\textwidth]{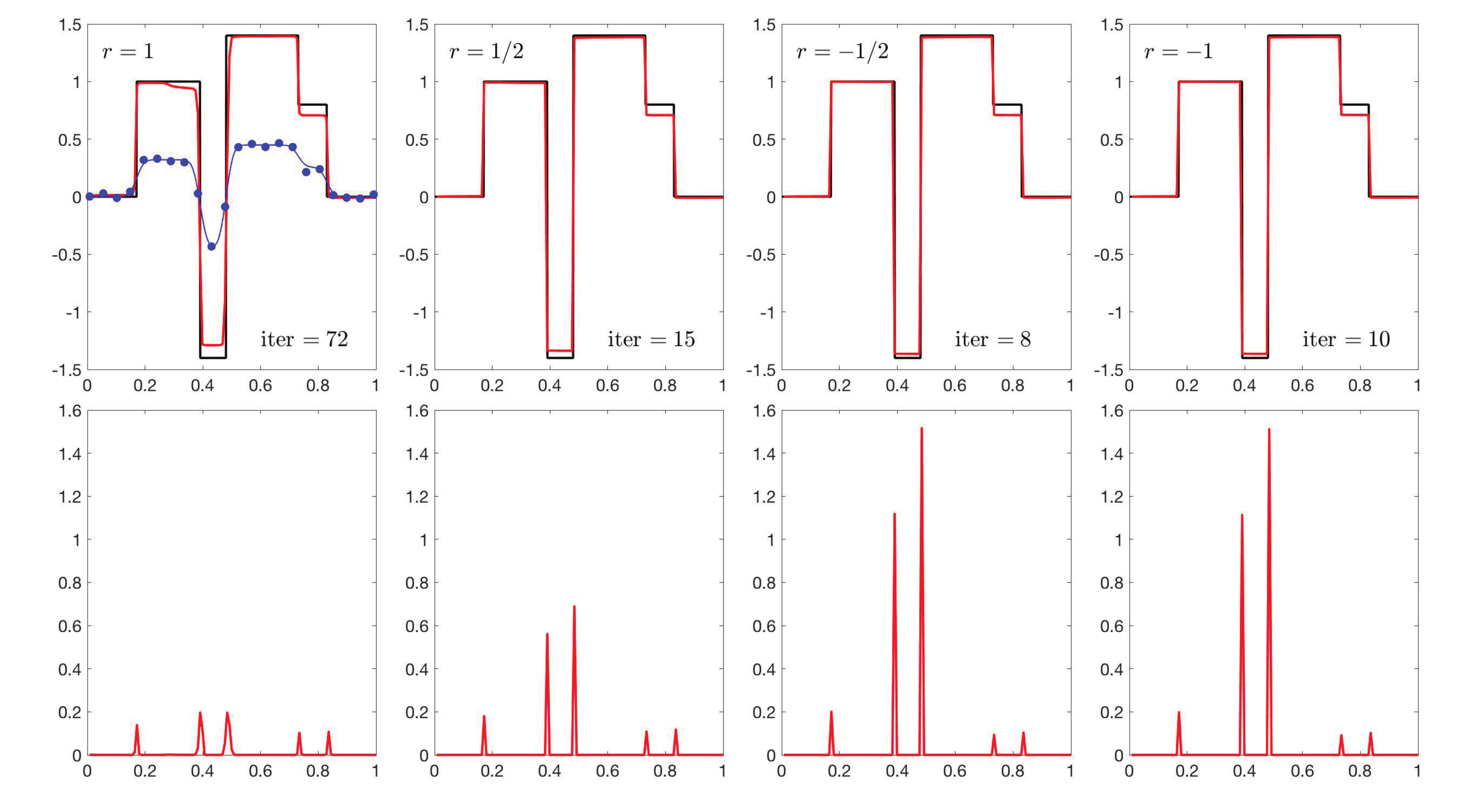}}
\caption{\label{fig:IAS results} The MAP estimates computed by the hybrid IAS algorithm for $z$ (upper row) in red, the black curves corresponding to the generative model, and for $\theta$ (lower row). The noisy data are shown as blue dots in the top left panel, the solid blue line indicating the noiseless convolution data.
 The left panels correspond to the gamma hyperprior $r=1$, or Phase I, and the result is the starting point for the Phase II iterations for $r<1$. The number of iterations for $r<1$ refers to the iteration rounds of Phase II. The parameter values are given in Table~\ref{tab:params}.}
\end{figure}

The MAP estimates computed by the IAS algorithm are shown in Figure~\ref{fig:IAS results}, where the number of iterations needed for satisfying the stopping criterion is indicated. Observe that the MAP estimate with $r=1$ is the starting point for all the hybrid models with $r<1$.

\subsection{Sampling}

We begin by applying the sampling algorithm with the gamma hyperprior ($r=1$) for the linear model (\ref{linear}), initiating the sample from the IAS-based MAP estimate. After some preliminary tests, the stepsize control parameter is set to $h = 0.05$, yielding consistently an acceptance rate close to 6.3\%.  Decreasing the stepsize increases the acceptance rate, e.g., $h = 0.02$ yields an acceptance rate of 33\%. The choice of the stepsize will be justified momentarily.

The relatively small stepsize $h$ implies that the draws in the sample are correlated, so to improve the sample quality, we retain only a subsample: In our test, we choose the computed sample size to be  10\,000\,000, and to decrease the dependency of the sample point, we keep only every 1000th point, reducing the effective sample size to $ N =10\,000$. The run time in a standard laptop is only 135 seconds, as the proposal density is pure white noise, and only one matrix-vector product for deciding on the acceptance is required. 

\begin{figure}[ht]
\centerline{\includegraphics[width = 0.5\textwidth]{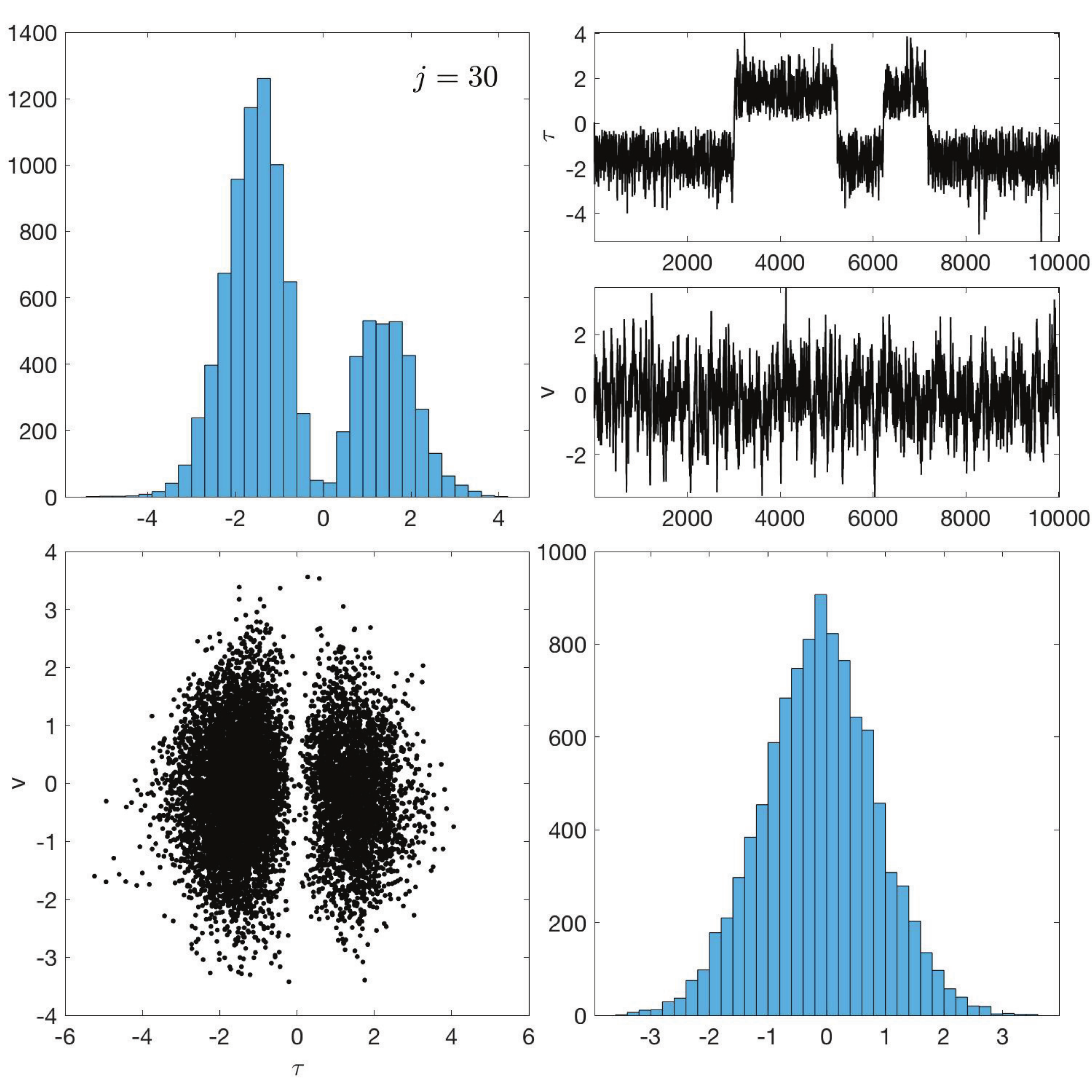}
\includegraphics[width = 0.5\textwidth]{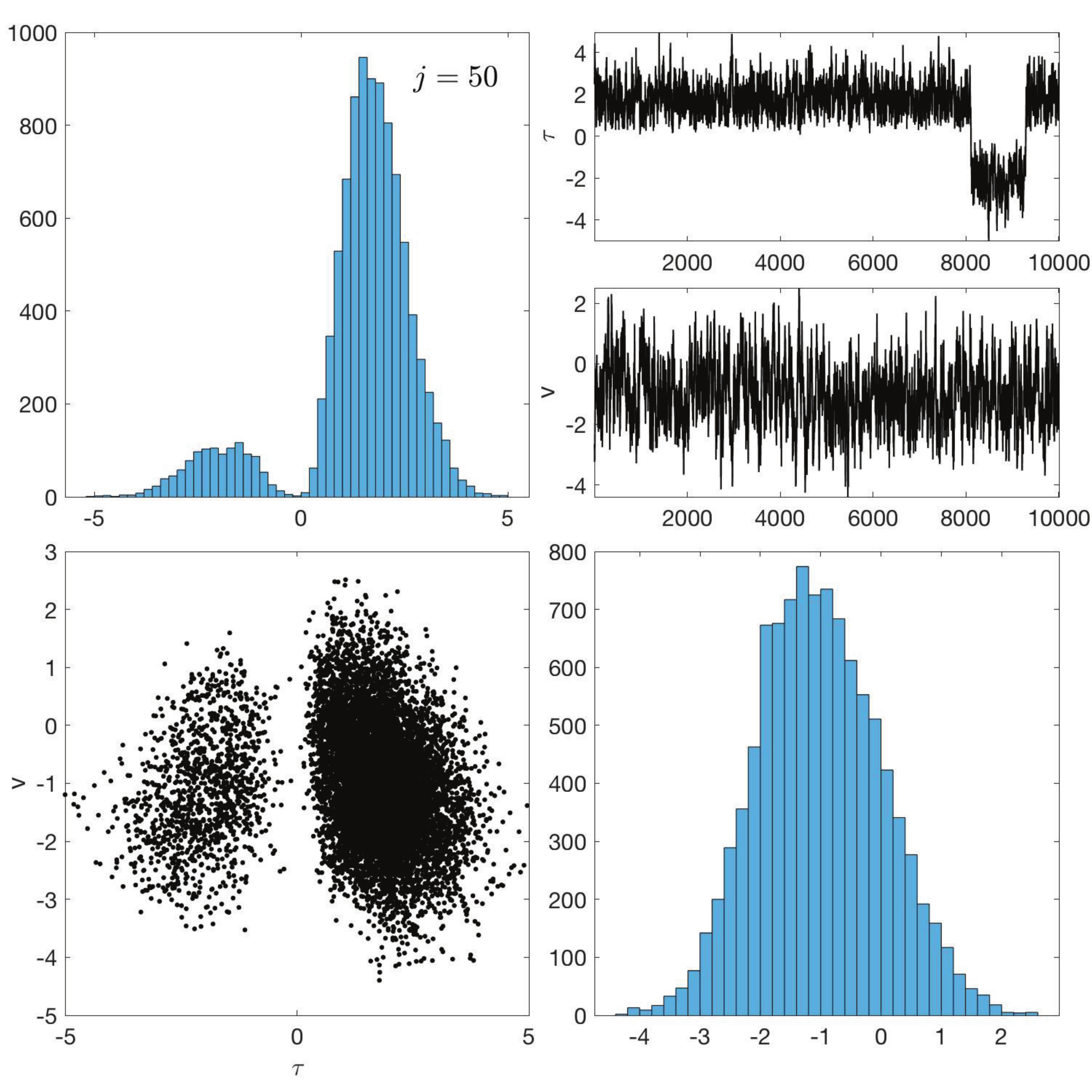}
}
\caption{\label{fig:Histograms gamma} Samples of two selected  pairs $(\tau_j,v_j)$ shown as scatter plots and histograms with hypermodel $r=1$. On the right, $j=30$,
corresponding to a position where the generative function is constant, and $j=50$ corresponding to a local minimum of the MAP estimate.}
\end{figure}

To analyze the mixing properties of the sampler, we select  two indices, $j_1 = 30$ and $j_2=50$, corresponding to values $t_{j_1} \approx 0.23$ and $t_{j_2} \approx 0.39$, the former corresponding to a position around which the generative function is constant, and the latter near a jump where the MAP estimate has a local minimum. Figure~\ref{fig:Histograms gamma}  shows the scatter matrices of the pair $(\tau_j,v_j)$ for $j=j_1$ and $j=j_2$, respectively, as well as the time traces of the sample.

We observe that the distribution of $\tau_j$ is bimodal, reflecting the fact that the physical scaled parameter $\lambda_j$ depends on $\tau_j^2$ which is insensitive to the sign of $\tau_j$. Furthermore, at $j=30$, the sample is centered near the coordinate origin, indicating that the sampler recognizes the point as belonging to the flat background. At $j=50$, the bimodal nature is visible, however, the values of $v_j$ are predominantly negative, indicating a presence of a negative jump, since
\[
 x_j = \vartheta^{1/2} \frac 1{\sqrt{2}}v_j |\tau_j|,\quad \theta_j = \frac 12 \vartheta \tau^2_j. 
\] 
However, the histogram of $v_j$ does not exclude the value $v_j=0$ corresponding to a background value $x_j=0$, indicating uncertainty in identifying the jump unambiguously.
Interestingly, if the stepsize $h$ is decreased to increase the acceptance rate, the sampler fails to identify the bimodal nature of the distribution and samples only from the mode $\tau_j>0$. This is not a serious issue, as the bimodality is simply a result of coordinate representation, however, we chose here select the strepsize so that the  feature becomes visible.

\begin{figure}[ht]
\centerline{\includegraphics[width = 0.5\textwidth]{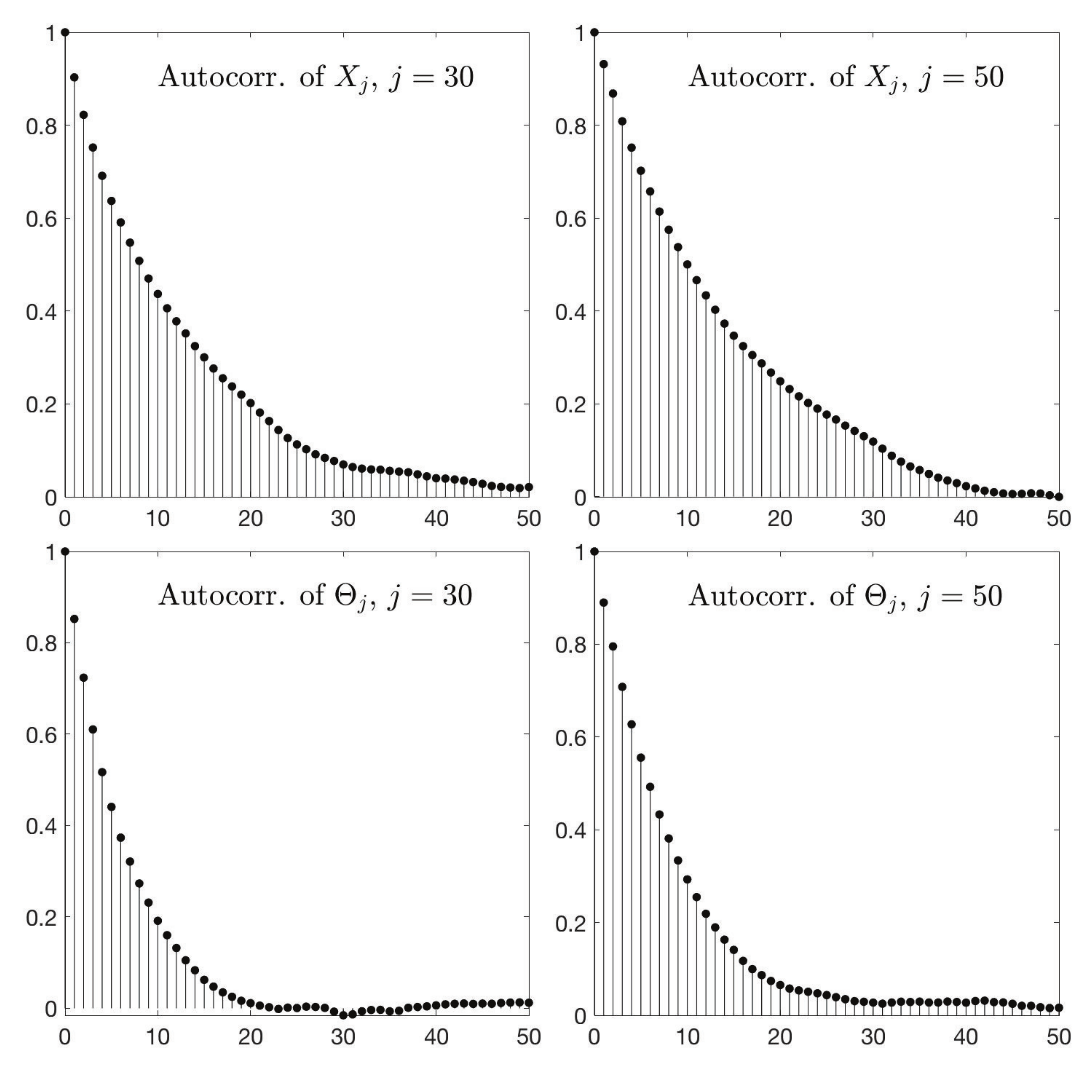}
\includegraphics[width = 0.5\textwidth]{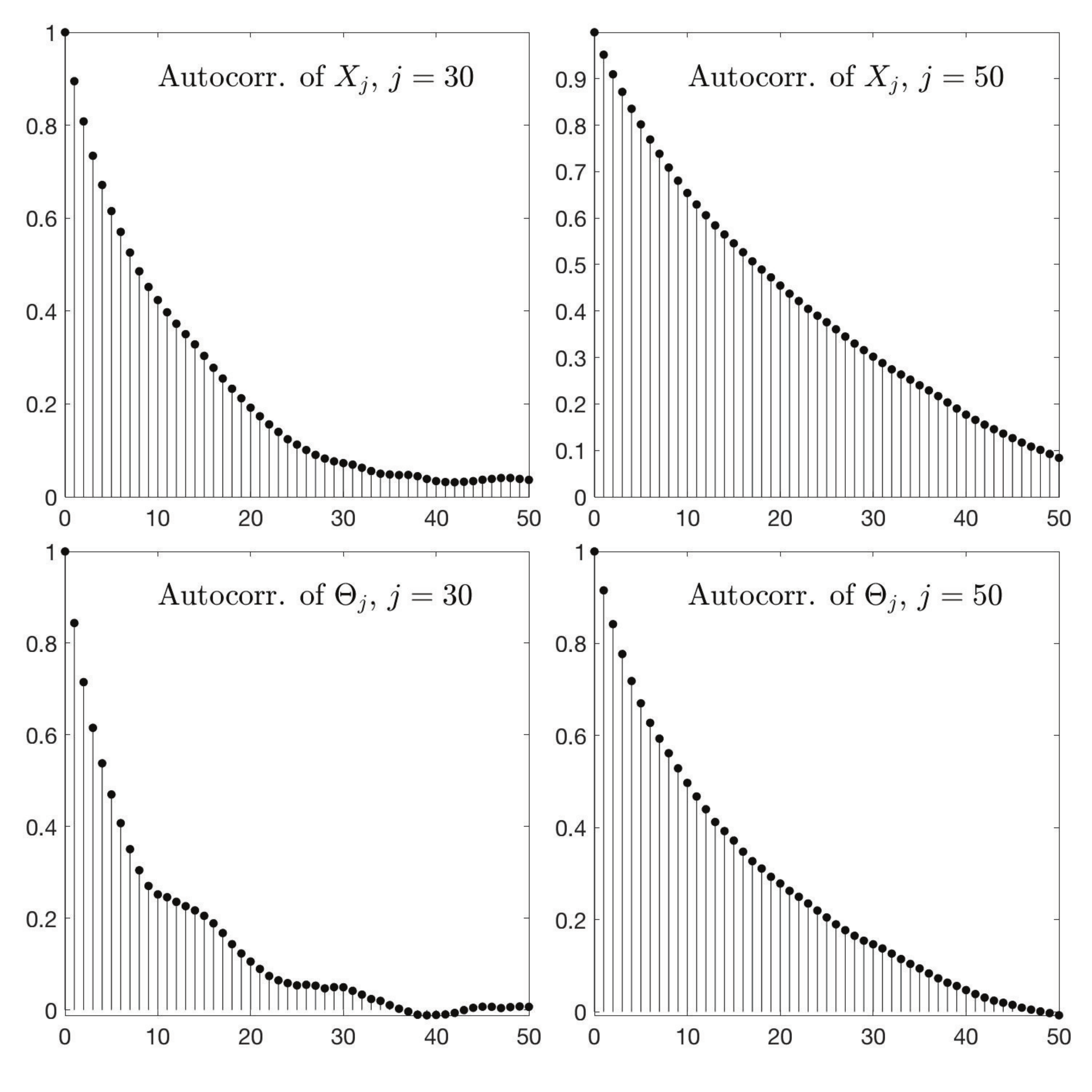}
}
\caption{\label{fig:autocorr gamma} Autocorrelation functions of the sample corresponding to the gamma hyperprior (two left columns) and to hypermodel with $r=1/2$ (two right columns). The top row shows the autocorrelation functions of the variables $X_j$ corresponding to a background value $j=30$ and to a jump value $j=50$, and the bottom row the autocorrelation function for the corresponding varibles $\Theta_j$.}
\end{figure}

Figure~\ref{fig:autocorr gamma} shows the autocorrelation functions of the retained sample of the components $(x_j,\theta_j)$, $j=30$ and $j=50$,
\[
 C(x_j)(\ell) = \frac 1{\|x_j\|}\sum_{k=1+\ell}^{N-\ell}(x_j^k - \overline x_j)(x_j^{k-\ell} - \overline x_j), \quad \ell =0,1,2,\ldots
\]
where
\[
 \overline x_j = \frac 1N\sum_{k=1}^N x_j^k, \quad \|x_j\| = \left(\sum_{k=1}^N \left(x_j^k\right)^2\right)^{1/2}.
\]  
The autocorrelation plots give a sense of the independence of the subsampled draws. There seems to be no significant difference between a background variable and one corresponding to a jump. 
  
\begin{figure}[ht]
\centerline{\includegraphics[width = 0.5\textwidth]{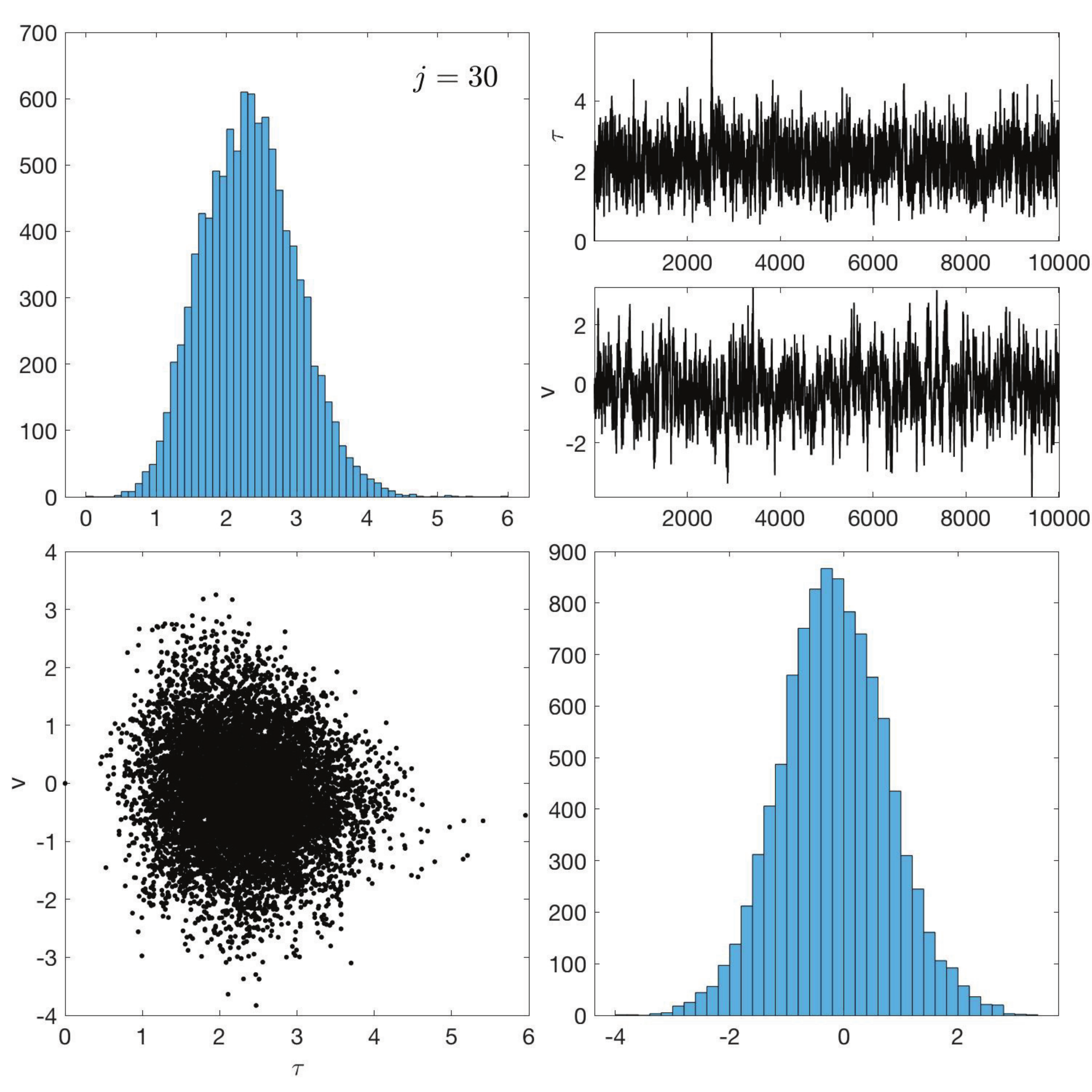}
\includegraphics[width = 0.5\textwidth]{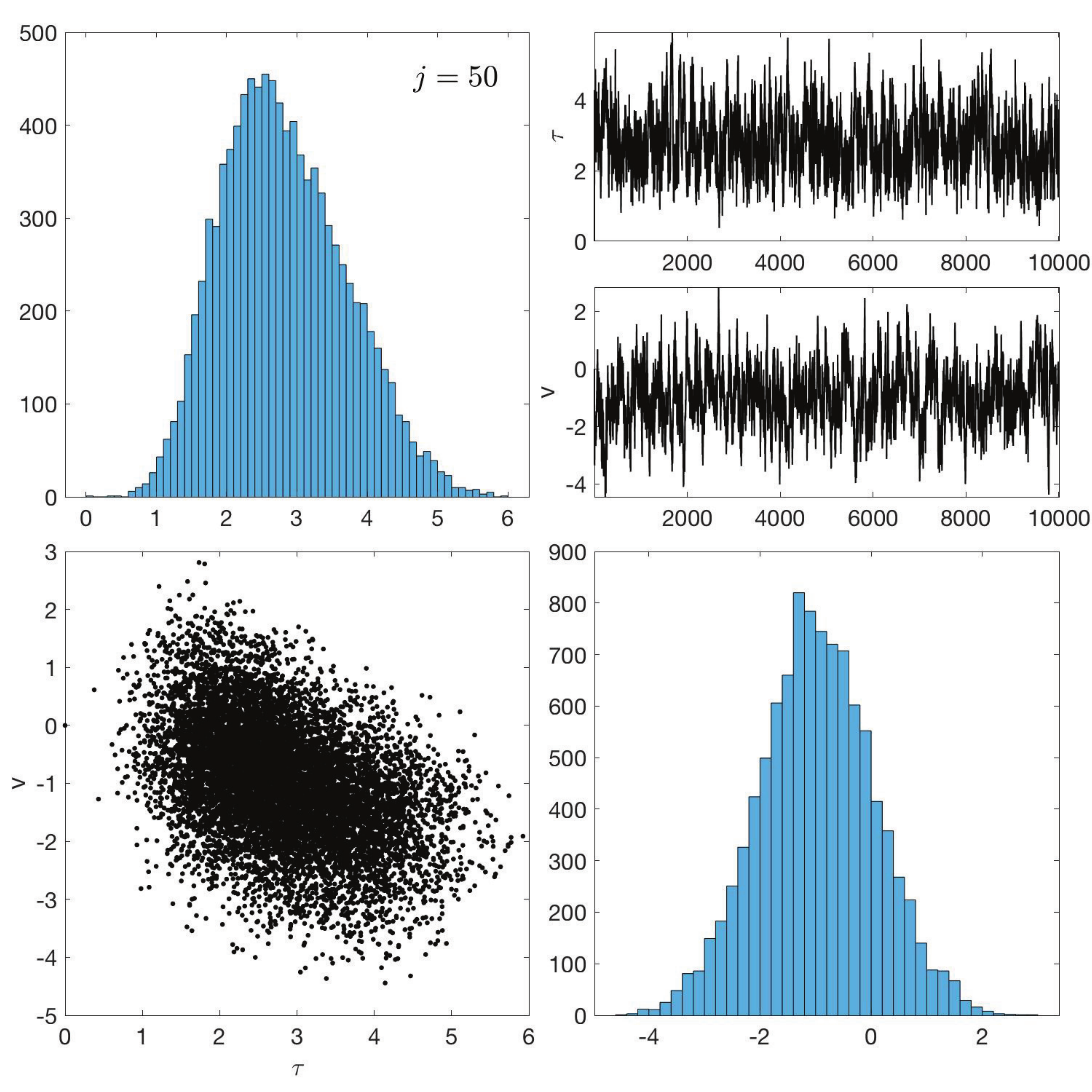}
}
\caption{\label{fig:Histograms half} Samples of two selected  pairs $(\tau_j,v_j)$ shown as scatter plots and histograms, hypermodel corresponding to $r=1/2$. On the right, $j=30$,
corresponding to a position where the generative function is constant, and $j=50$ corresponding to a local minimum of the MAP estimate.}
\end{figure}

Before analyzing the sample further, we run similar pCN sampling using the other hypermodels. It turns out that as $r$ decreases, the sampling becomes more challenging. We start with $r=1/2$. Using the same stepsize $h=0.05$ as in the case $r=1$ the acceptance rate falls to 4.8\%, so we decrease the stepsize slighltly to $h=0.03$, yielding an acceptance rate of 16.1\%. Generating a sample of size 10\,000\,000, retaining again every 1000th sample point takes 131 seconds in a standard laptop. In this case, even with a larger stepsize, the sampler is unable to detect the coordinate bimodality and samples only from the mode where $\tau_j>0$, as shown by the scatter plots in Figure~\ref{fig:Histograms half}.
 
The level of independence of the draws can be inferred from the autocorrelation functions shown in Figure~\ref{fig:autocorr gamma}. We observe that the correlation level is not significantly different from that corresponding to the gamma hypermodel.

Consider now the hypermodels with negative $r$, yielding a highly non-linear Gibbs energy functional with a strong sparsity promotion in the MAP estimation problems. It turns out that these models pose a challenge for sampling as well. We start with $r = -1/2$. Using the same stepsize as in the case $r=1/2$, $h=0.03$ the acceptance rate is as low as 0.12\%, so it is reasonable to decrease the stepsize. After extensive testing, we found that $h =0.008$ yields a reasonable acceptance rate of 6\%, and $h=0.005$ gives a 12\% acceptance.  Since the former figure is close to the values in the previous experiments, we set $h = 0.008$. To make the results comparable to the previous ones, we keep the same sample size of 10\,000\,000, retaining every 1000th sample point. The computing times are close to the ones reported above, of the order of 110 seconds.

\begin{figure}[ht]
\centerline{\includegraphics[width = 0.5\textwidth]{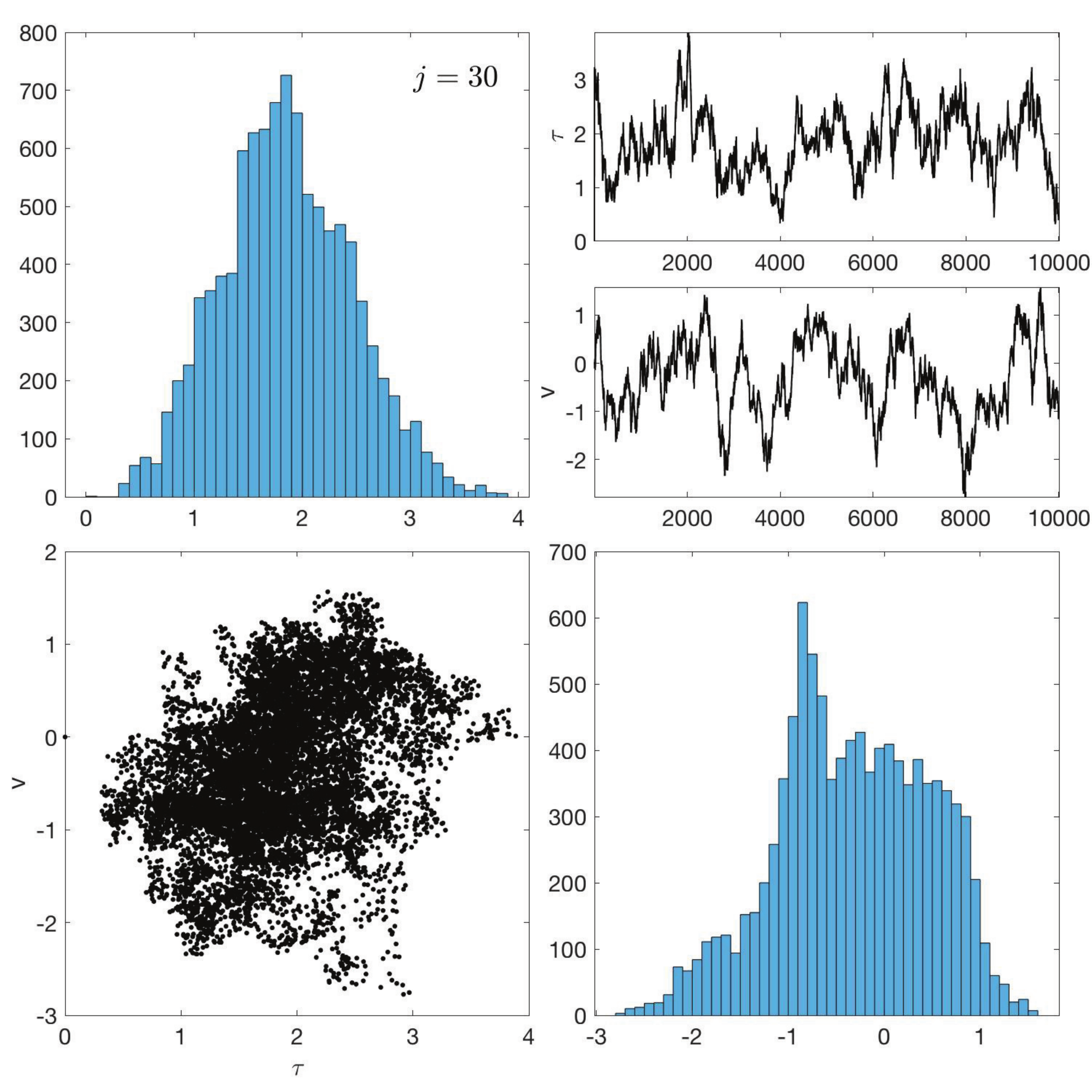}
\includegraphics[width = 0.5\textwidth]{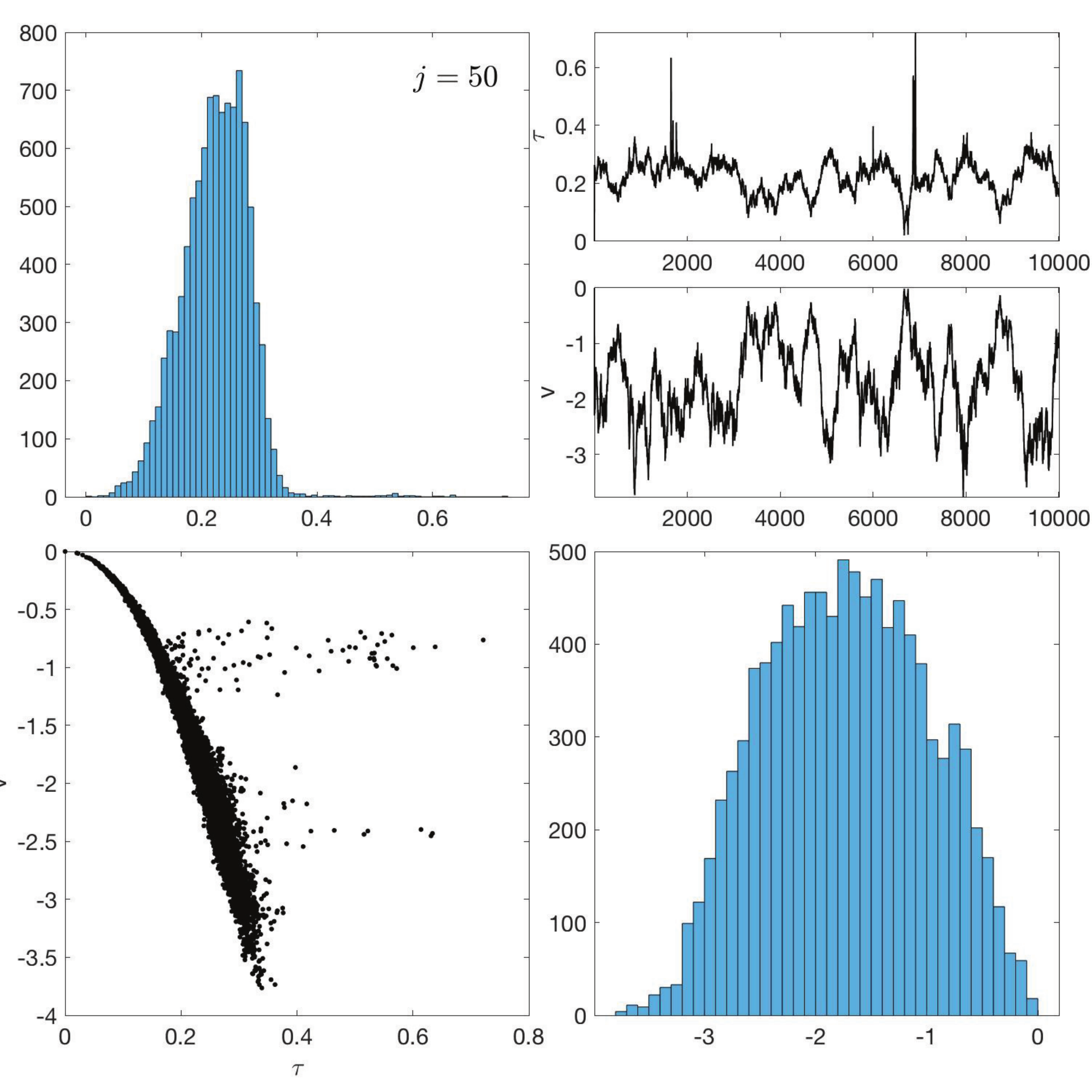}
}
\caption{\label{fig:Histograms minus half} Samples of two selected  pairs $(\tau_j,v_j)$ shown as scatter plots and histograms with hypermodel $r=-1/2$. On the right, $j=30$,
corresponding to a position where the generative function is constant, and $j=50$ corresponding to a local minimum of the MAP estimate.}
\end{figure}  

Figure~\ref{fig:Histograms minus half} shows the scatter plots and histograms of the parameters $(\tau_j,v_j)$. We observe that at $j=50$, the points are concentrated near a parabolic curve, which is to be expected, as the variables are related to each other through the formula
\[
 \xi_j = 2\frac{v_j}{\tau_j^2},
\]
and at $j=50$, we expect that the likelihood favors combinations of $(\tau_j,v_j)$  yielding a negative value for  $\xi_j$. 

\begin{figure}[ht]
\centerline{\includegraphics[width = 0.5\textwidth]{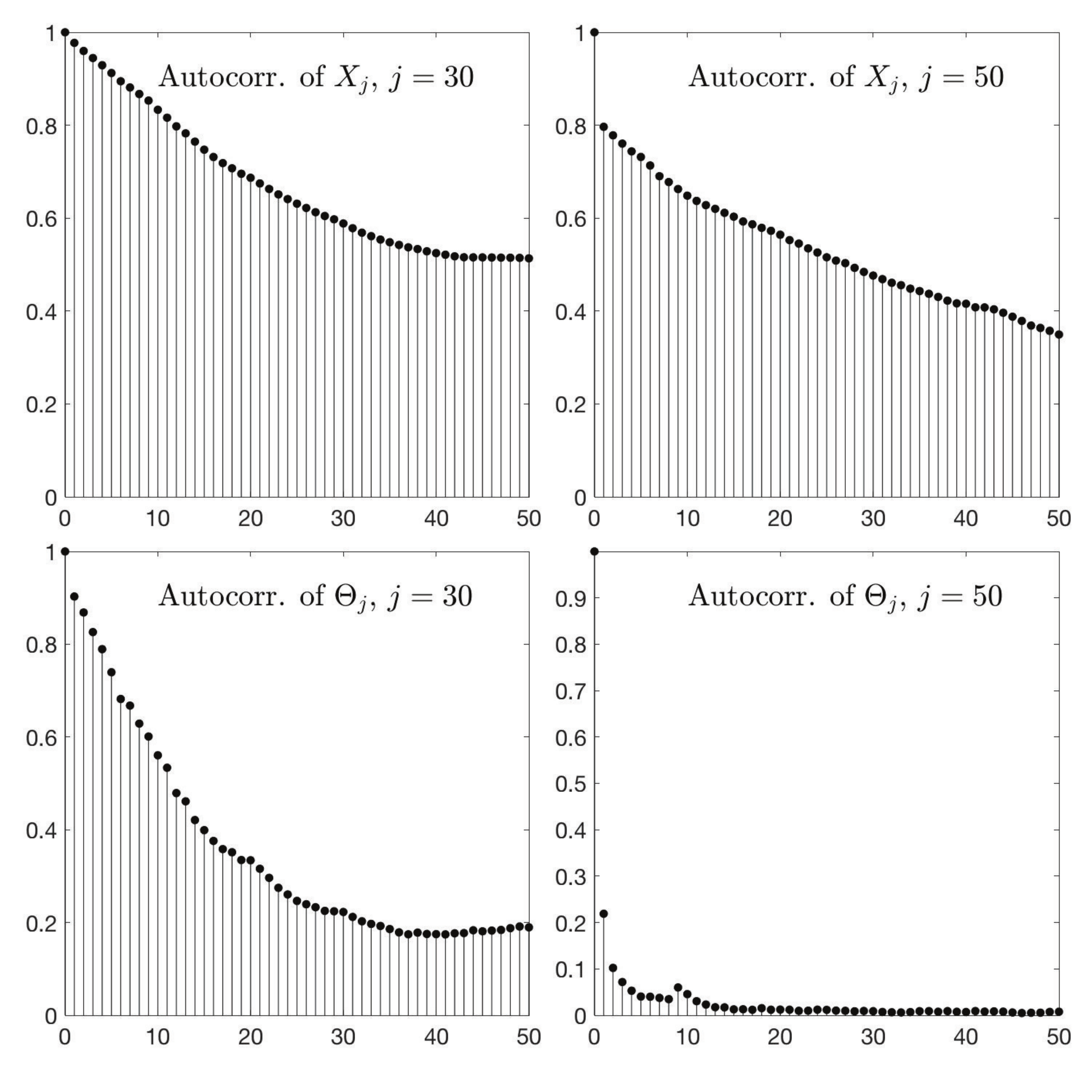}
}
\caption{\label{fig:autocorr negative} Autocorrelation functions of the sample corresponding to the hyperprior with $r=-1/2$. The top row shows the autocorrelation functions of the variables $X_j$,  $j=30$, corresponding to a background value and to a jump value, $j=50$; the bottom row shows the autocorrelation functions for the corresponding variables $\Theta_j$.}
\end{figure}

The sample histories of $\tau_j$ and $v_j$ are not indicating a good mixing, however, from the practical point of view, the important question is how well the sample mixes the values $(x_j,\theta_j)$. To explore this question, we plot again the autocorrelation functions of the selected components in Figure~\ref{fig:autocorr negative}. Interestingly, while the autocorrelation functions of the $x$-variables are decreasing rather slowly, the $\theta$-variable autocorrelations indicate rather good mixing.

\begin{figure}[ht]
\centerline{\includegraphics[width = 0.5\textwidth]{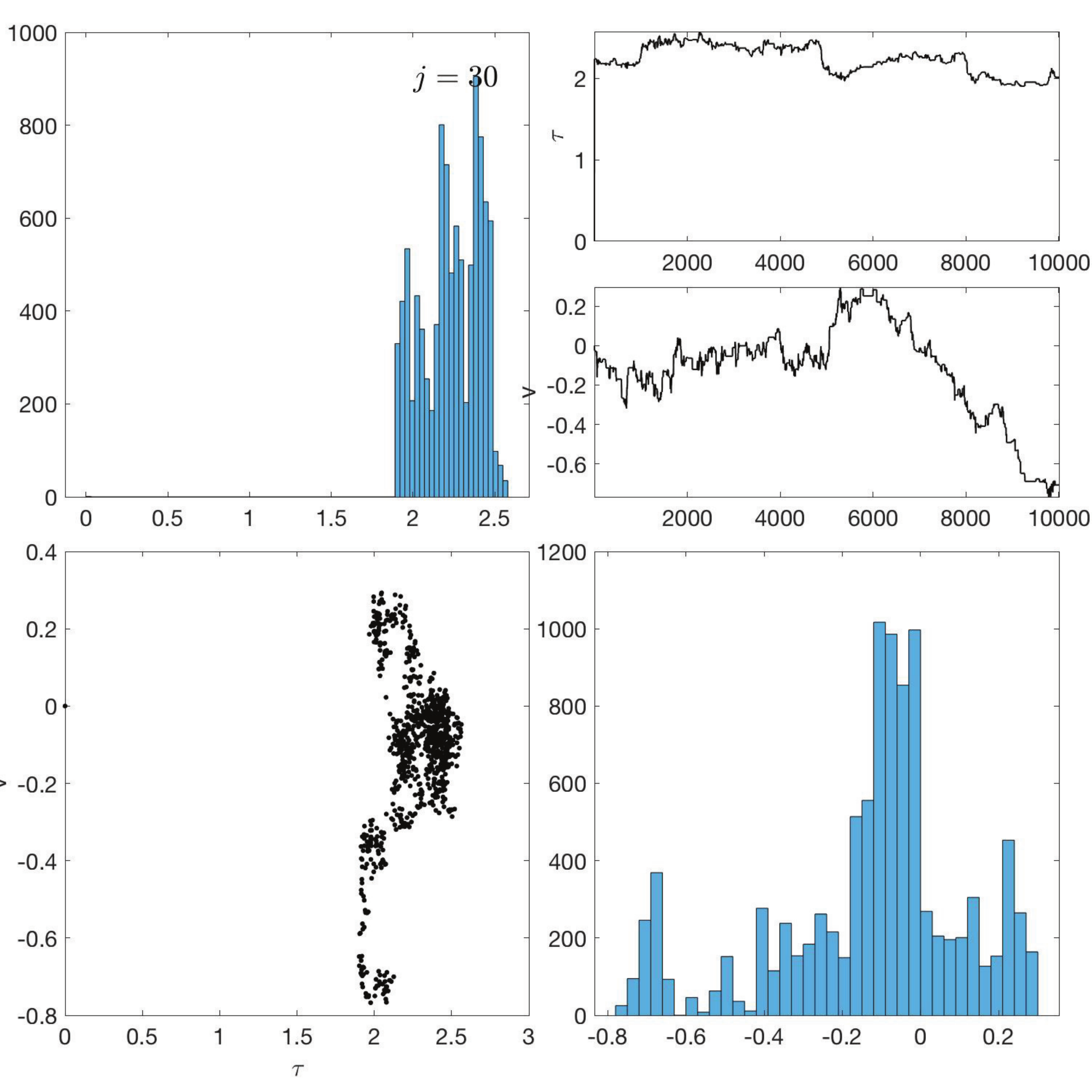}
\includegraphics[width = 0.5\textwidth]{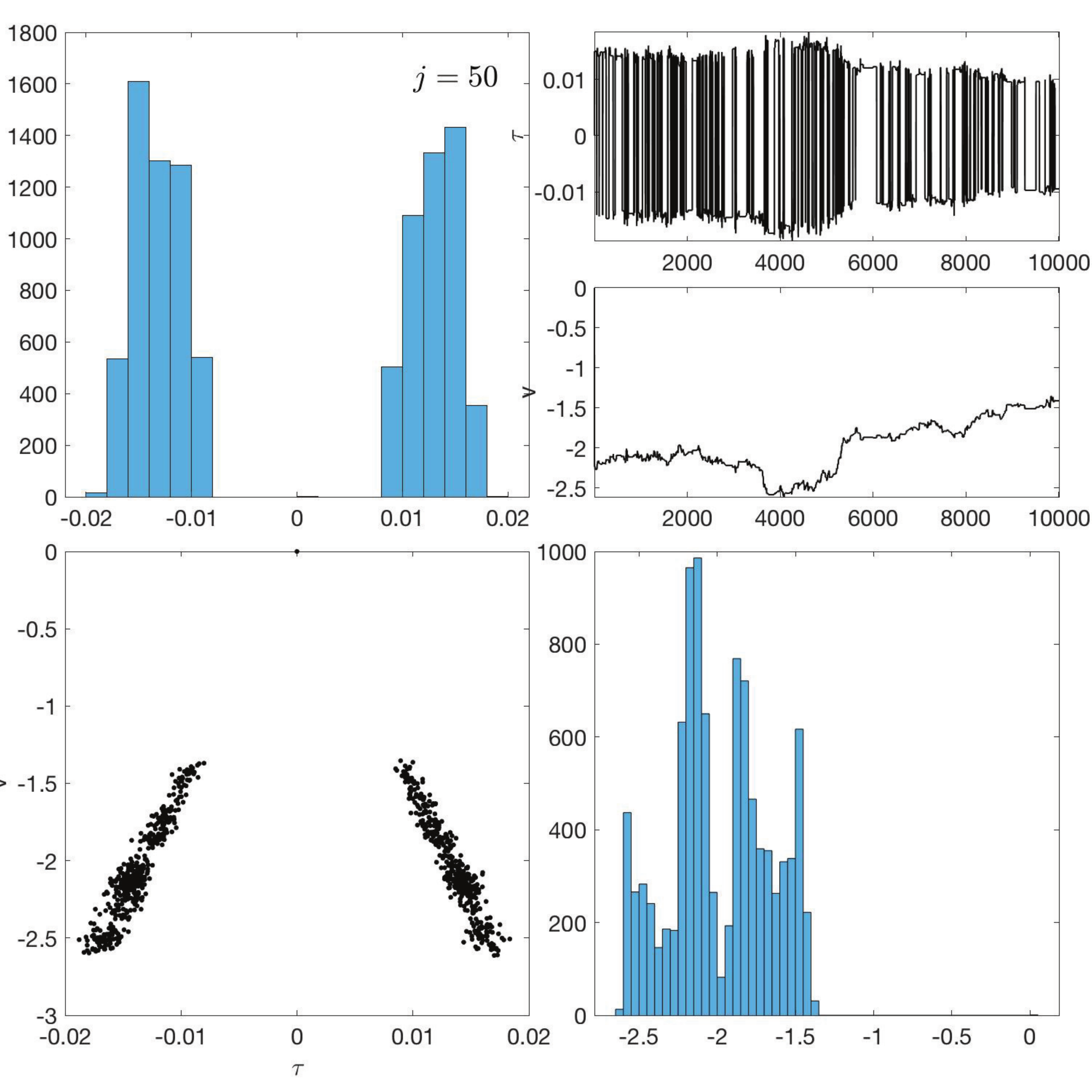}
}
\caption{\label{fig:Histograms inv gamma} Samples of two selected  pairs $(\tau_j,v_j)$ shown as scatter plots and histograms with hypermodel $r=-1$. On the right, $j=30$,
corresponding to a position where the generative function is constant, and $j=50$ corresponding to a local minimum of the MAP estimate.}
\end{figure}

Finally, we consider the inverse gamma hyperprior, $r=-1$. Numerical tests indicate that the proposed sampler struggles to find a reasonably well mixed sample. To demonstrate this, we 
select first the step size to be $h=0.02$, which leads to an acceptance rate as low as 0.001\%, that is, one proposal of every 100\,000 is accepted in the average. We then generate a sample of size $10^8$, keeping only every 10\,000th point, i.e., we have a sample of size 10\,000 in which approximately 90\% of points are repeated values corresponding to rejections. The computing time of this sample is less than 18 minutes. The $(\tau_j,v_j)$ scatter plots for the two selected values of $j$ are shown in Figure~\ref{fig:Histograms inv gamma}.

The time traces of the samples reveal that the sample is of low quality, and reliable conclusions could hardly be based on this sample. Numerical experiments show that decreasing the step size does not improve significantly the sample quality.  The plots show, however, some features that may help designing a better sampler.  First, the scatter plot for $j=50$ shows again that the coordinate transformation leads to a bimodal distribution, which is hard to reproduce with a smaller step size. Second, we observe that the sampler proposes points that lie along lines passing through the origin. Recalling that for $r=-1$, the coordinate transformation implies that
\[
\xi_j = \sqrt{2} \frac{v_j}{|\tau_j|},
\]
the absolute value of the slope of the line determines the proposed value for $\xi_j$. It turns out that for $j=50$, the sample history of $\xi_{50}$ is quite satisfactory, with a relatively short correlation length (not shown here.)

To improve the mixing in the case $r=-1$, we propose a modification of the pCN algorithm based on a reparametrization of the pairs $(\tau_j,v_j)$. Consider the problem of generating a standard normal distribution in the plane $\R^2$.  Let $X^{j-1}\sim {\mathcal N}(0,\mI_2)$, and denote $R^{j-1} = \|X^{j-1}\|$. We have
\[
 R^{j-1} \sim{\rm Rayleigh}(1).
\]
Let $W\sim {\mathcal N}(0,\mI_2)$, and assume that $W$ is independent of $X^{j-1}$. For any $k>0$, we define
\begin{equation}\label{radial prop}
 R^j = \left( (1-k^2)\left(R^{j-1}\right)^2 + 2 k \sqrt{1-k^2} R^{j-1} W_1 + k^2 \|W\|^2 \right)^{1/2}.
\end{equation}
We claim that
\[
 R^j   \sim{\rm Rayleigh}(1).
\] 
This is a direct consequence of the fact that
\[
 X^j = \sqrt{1-k^2} X^{j-1} + k W \sim {\mathcal N}(0,\mI_2).
\]
Without loss of generality, we may assume that the coordinates are chosen so that $X^{j-1}_2 = 0$. Then
\[
 R^j = \|X^j\|,
\]
follows a Rayleigh distribution.

Let $\Phi^{j-1}$ be the phase angle of $X^{j-1}$,
\[
 \Phi^{j-1} = {\rm atan}(X_2^{j-1},X_1^{j-1}) \sim {\rm Uniform}(\R/2\pi).
\] 
For an arbitrary $h>0$, define
\begin{equation}\label{angular prop}
 \Phi^j = \Phi^{j-1} + h \Omega, \quad \Omega\sim{\mathcal N}(0,1).
\end{equation}
Then, $\Phi^j \sim {\rm Uniform}(\R/2\pi)$. Based on these observations, we introduce the following two-phase proposal:

{\em Given $h>0$, $k>$, and the current point $x^{j-1} = (\tau^{j-1},v^{j-1})\in\R^2$,}
\begin{enumerate}
\item {\em Set $r^{j-1} = \|x^{j-1}\|$, and}
\[
  r^j = \left( (1-k^2)\left(r^{j-1}\right)^2 + 2 k \sqrt{1-k^2} r^{j-1} w_1 + k^2 \|w\|^2 \right)^{1/2}, \quad w\sim {\mathcal N}(0,\mI_2);
\]
\item {\em Set    $\varphi^{j-1} = {\rm atan}(v^{j-1},\tau^{j-1})$, and}
\[
 \varphi^j = \varphi^{j-1} + h \omega, \quad \omega\sim{\mathcal N}(0,1).
\]
\item{\em  Define}
\[
 x^j = (\tau^j,v^j) = (r^j\cos\varphi^j,r^j\sin\varphi^j).
\]  
\end{enumerate}

In this manner, the new variable $X^j$ with realization $x^j$ follows the same Gaussian distribution as $X^{j-1}$ as in the standard pCN proposal, but we control separately the stepsize in the radial and in the angular direction. The two free variables $k$ and $h$ add complexity to the tuning process, but may lead to a chain with better mixing properties. 

\begin{figure}[ht]
\centerline{\includegraphics[width = 0.5\textwidth]{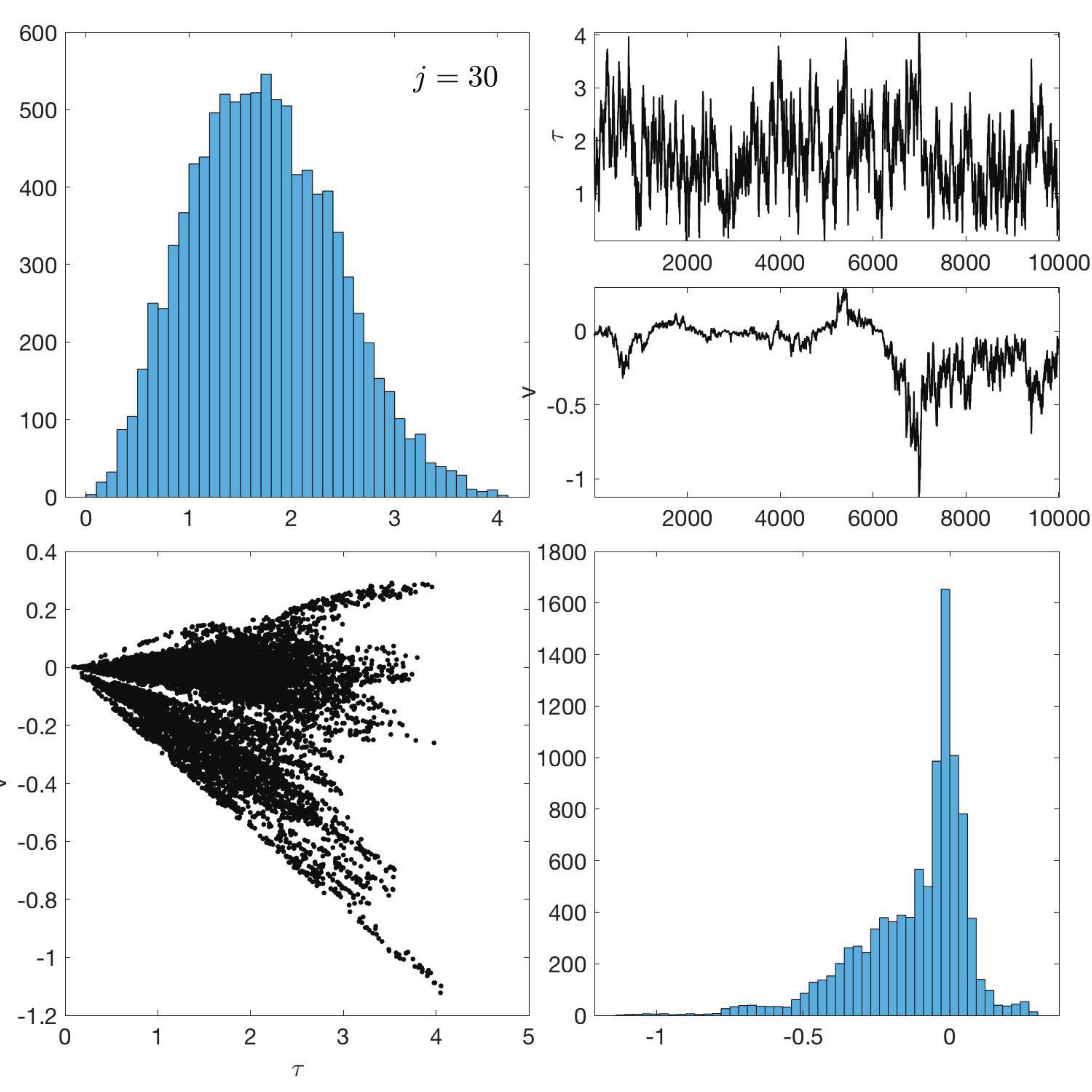}
\includegraphics[width = 0.5\textwidth]{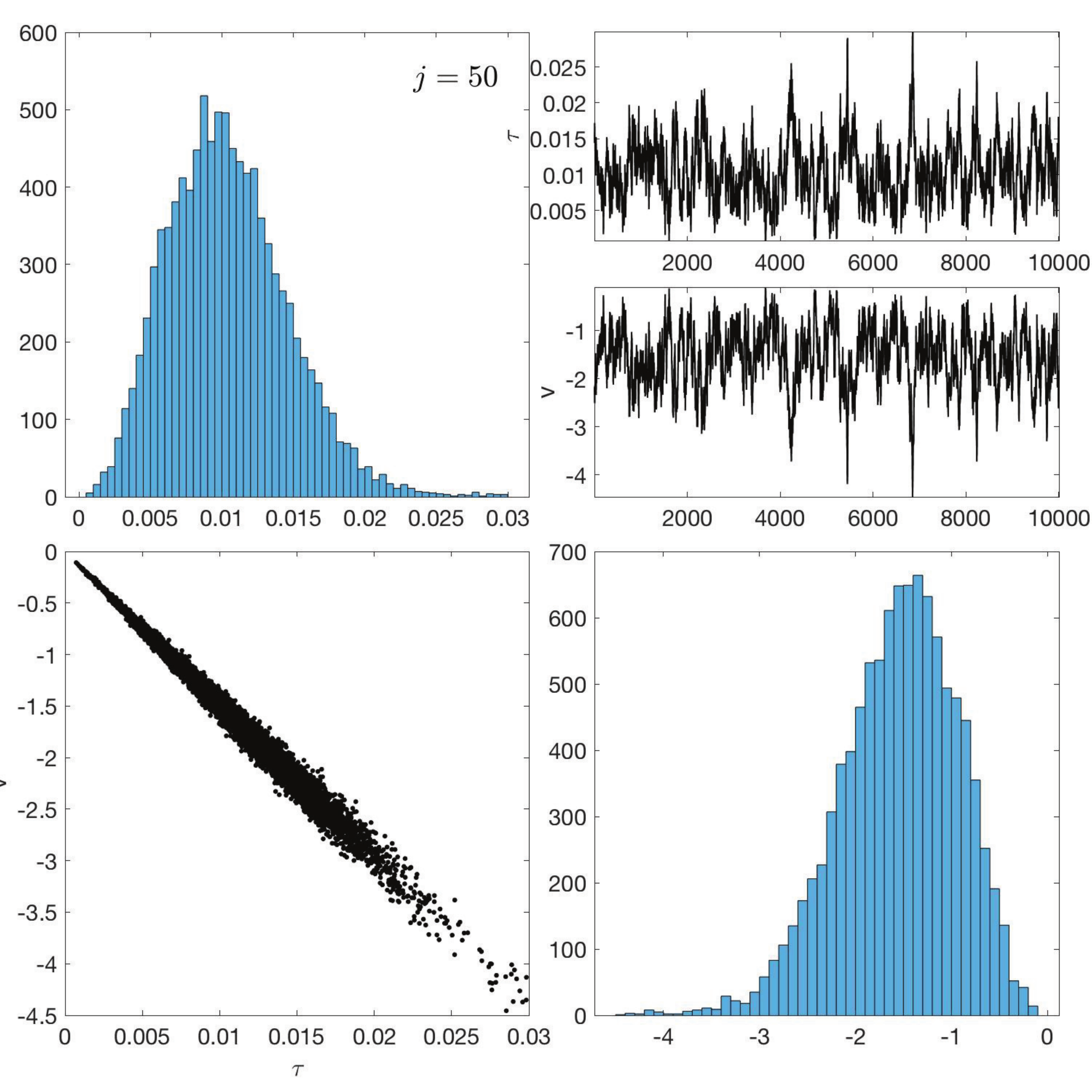}
}
\caption{\label{fig:Histograms radial} Samples of two selected  pairs $(\tau_j,v_j)$ shown as scatter plots and histograms with hypermodel $r=-1$ using the modified pCN algorithm. On the right, $j=30$, corresponding to a position where the generative function is constant, and $j=50$ corresponding to a local minimum of the MAP estimate.}
\end{figure}

To demonstrate the viability of the proposed algorithm, we run the sampler with parameter values $h = 0.001$ and $k = 0.05$, yielding an acceptance rate of 1.5\%. We compute a sample of size 5\,000\,000, retaining every 500th realization. The computing time is slightly longer than with plain pCN, requiring 91 seconds on a standard laptop. The scatterplots corresponding to this sample are shown in Figure~\ref{fig:Histograms radial}.

\begin{figure}[ht]
\centerline{\includegraphics[width = 0.5\textwidth]{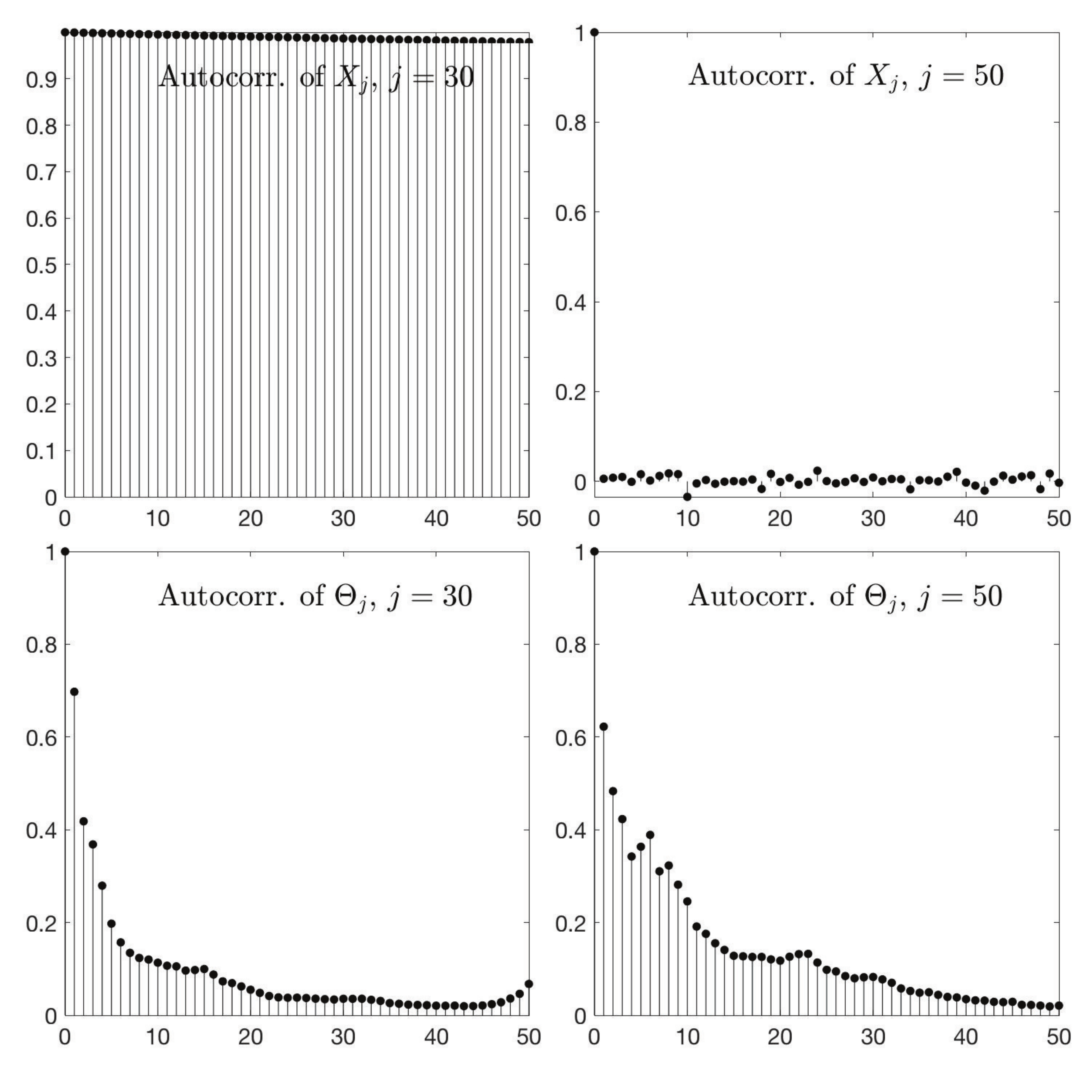}
}
\caption{\label{fig:autocorr radial} Autocorrelation functions of the sample corresponding to the hyperprior with $r=-1$  corresponding to the modified pCN sampler. The top row shows the autocorrelation functions of the variables $X_j$ corresponding to a background value $j=30$ and to a jump value $j=50$, and the bottom row for the corresponding varibles $\Theta_j$.}
\end{figure}

To assess the quality of the sample, we compute again the autocorrelation functions for the selected variables, see Figure~\ref{fig:autocorr radial}. We observe that while the autocorrelation indicates poor quality of the sample of $X_{50}$, the autocorrelation of $\Theta_{50}$ decreases relatively rapidly. The conclusion therefore is that if the sample-based estimate of $\Theta_{50}$ is small, we may claim with high certainty that $X_{50}$ is small, too.
 
\begin{figure}
\centerline{\includegraphics[width=12cm]{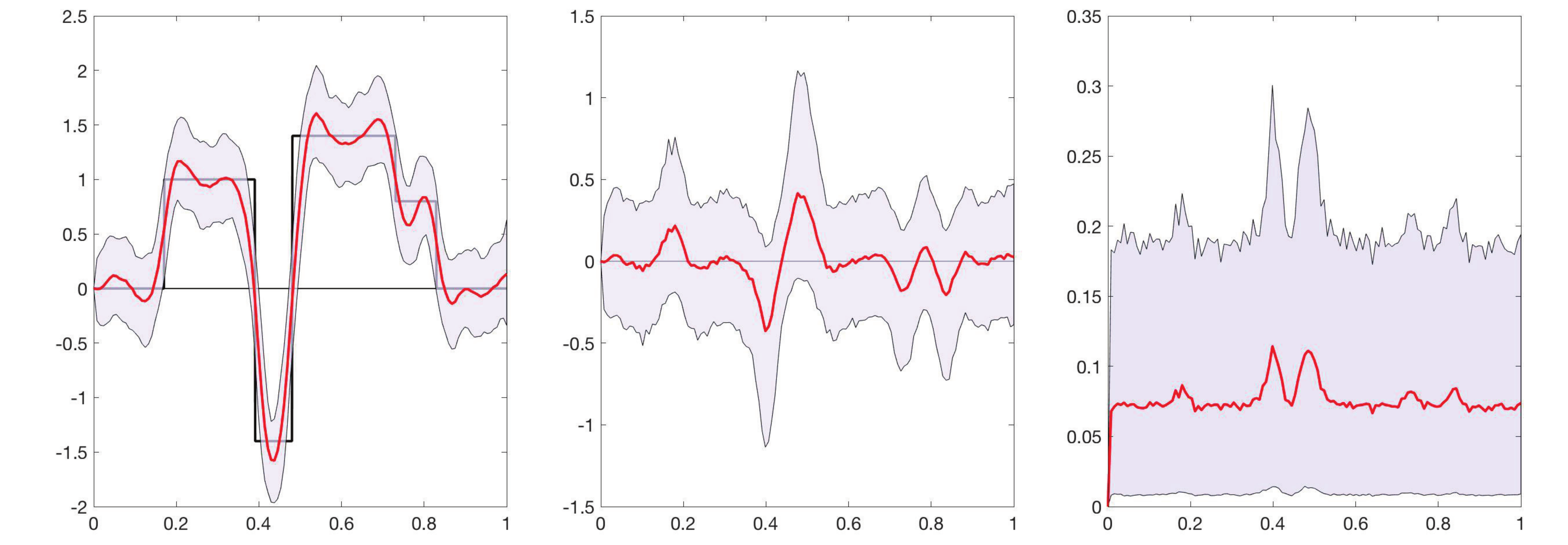}}
\centerline{
\includegraphics[width=12cm]{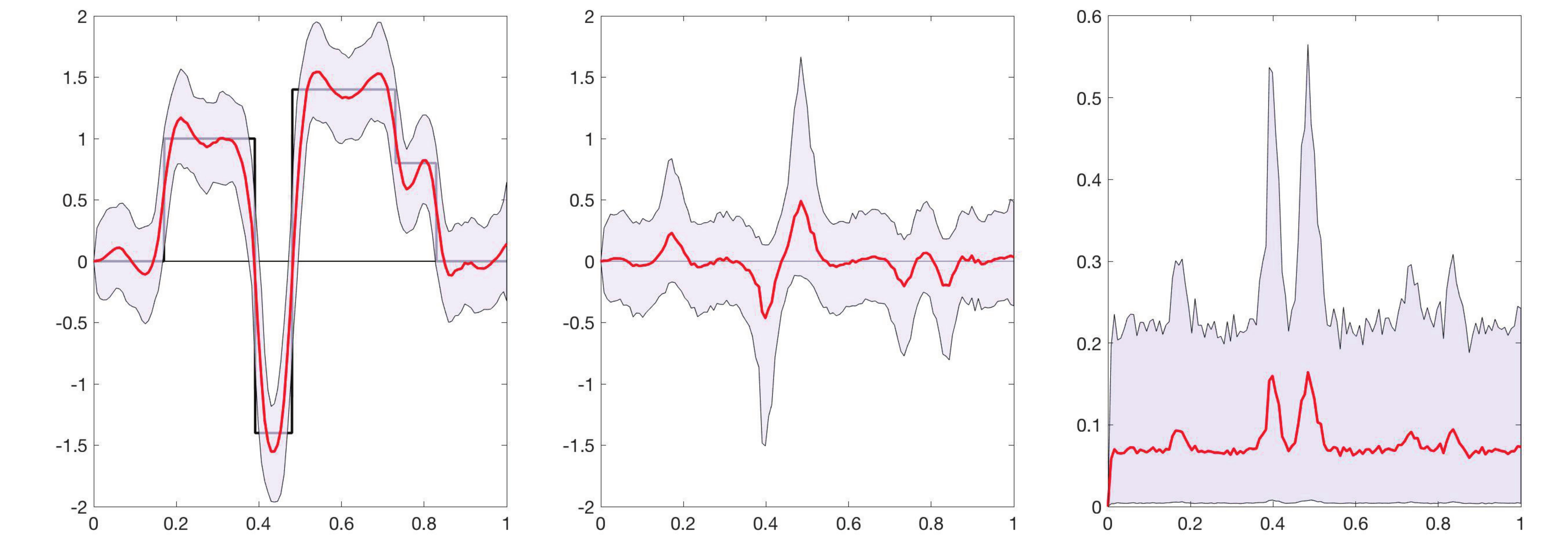}
}
\centerline{
\includegraphics[width=12cm]{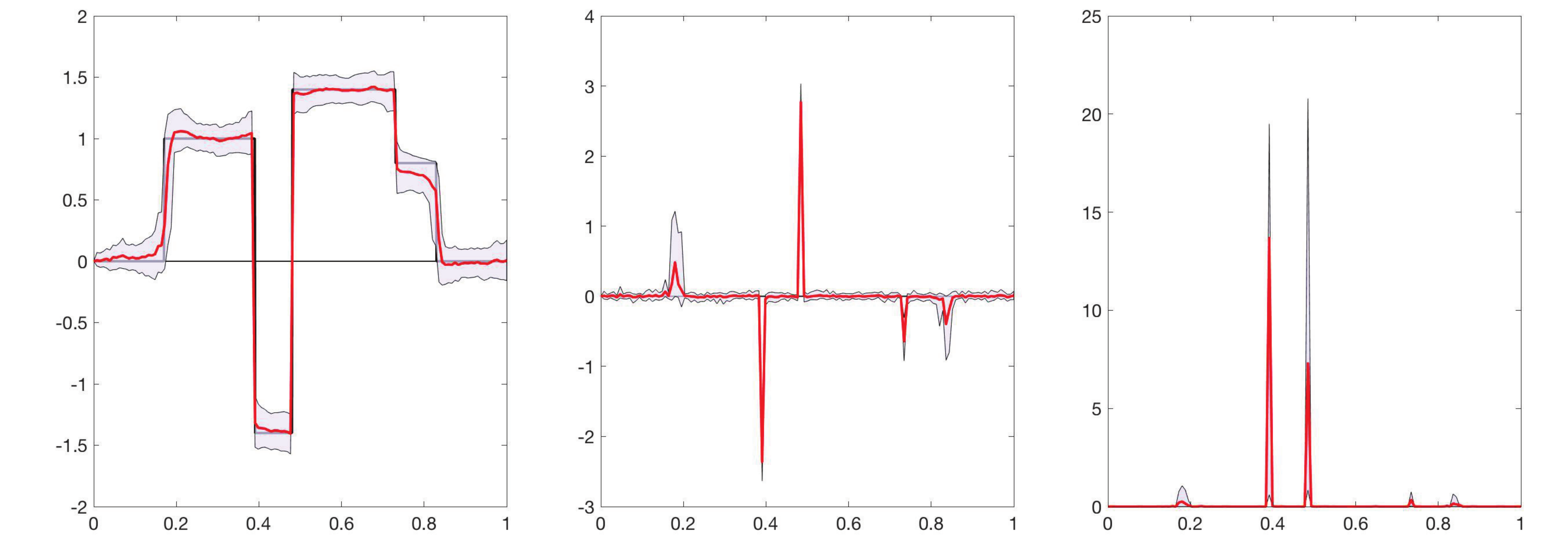}
}
\centerline{
\includegraphics[width=12cm]{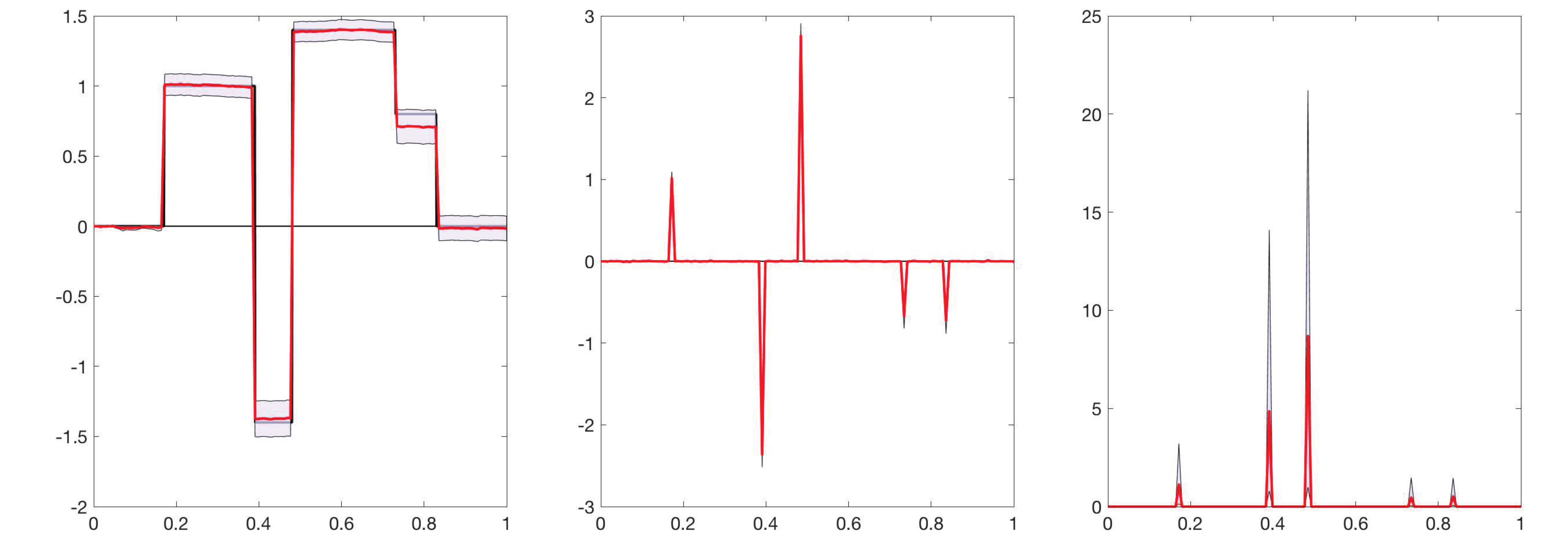}
}
\caption{\label{fig:envelopes} Posterior means (red curve) and the 90\% credible envelopes for $z$ (left), $x$ (middle) and $\theta$ (right) with hypermodels $r=1$ (top), $r=1/2$, $r=-1/2$ and $r = -1$ (bottom). The inverse gamma sample is based on the modified pCN algorithm.}
\end{figure}

To summarize the results of the samplers, in Figure~\ref{fig:envelopes} we plot the estimated posterior means and 90\% credible envelopes of the variables $z$, $x$ and $\theta$ corresponding to the four hypermodels.

The results show that while the MAP estimates for all models can be considered sparse solutions, the posterior means with $r=1/2$ and $r=1$ in particular, do not reflect the sparsity promoting nature of the prior models. This suggests that the posterior mean is not a particularly good representative summary of the posterior distribution, even if the individual samples would reflect the sparsity prior. 

\begin{figure}
\centerline{\includegraphics[width=\textwidth]{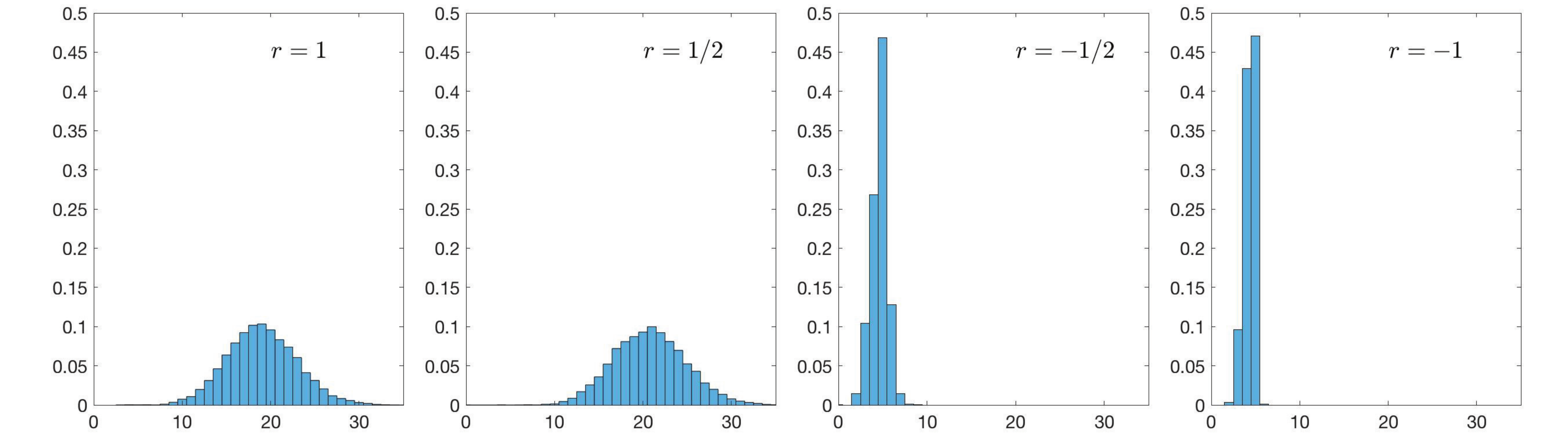}}
\caption{\label{fig:sparsity} Histograms of the number of components in the vectors $\theta^j$ that are above the threshold value $\delta$ (\ref{threshold}), thus indicating the level of compressibility of the sample vectors.  Observe that the number of non-zero increments in the generative model is 5, corresponding to the maximum for  $r=-1/2$ and $r=-1$.}
\end{figure}

The smoothness of the posterior mean raises the question to what extent  profiles sampled from the posterior density are compressible. To further investigate this issue, let $\theta^j$ denote the $j$th sample vector of the variance parameter, and define $\delta$-compressibility by the formula
\[
 \|\theta^{j}\|_{0,\delta} = {\rm card}\big\{1\leq k\leq n\mid \theta^{j}_k >\delta\big\},
\]
where $\delta>0$ is a given threshold. We set the threshold to correspond to one standard deviation above the mean  of the gamma distribution,
\begin{equation}\label{threshold}
 \delta = \beta_1 \vartheta + \sqrt{\beta_1}\vartheta.
\end{equation}
Figure~\ref{fig:sparsity} shows the histograms of the number of components in the vectors $\theta^j$ exceeding the threshold $\delta$. The samples based on models $r=1$ and $r=1/2$ identify significantly more increments above the threshold than the cardinality of the support of the generative model, the maximum of the histograms being around twenty, while for both $r=-1/2$ and $r=-1$, the maximum of the histogram is at 5, which coincides with the number of non-zero increments in the generative model.

\section{Discussion}

The difficulties of sampling from posterior densities corresponding to hierarchical Bayesian models are twofold: The curse of dimensionality makes efficient sampling hard, and the strong correlation between the parameters at different level adds an extra bottleneck for samplers. In this article, the former problem has been addressed by transforming the problem so that the preconditioned Crank-Nicholson sampling scheme can be applied. As the examples with strongly non-convex energy functionals show, the transformation does not completely remove the second problem, however, the transformation gives enough insight to allow further developments of the sampling strategy so that at least in the case of inverse gamma hyperprior model, a relatively efficient algorithm can be found. A natural question to be addressed in the future is if the proposed algorithm can be generalized or hyperprior models more general than the inverse gamma distribution. 

The sampling analysis for sparsity promoting hypermodels reveals that the concept of sparsity promotion is more complex than the analysis of the maximum a posterior estimates reveals: The MAP estimate may identify the correct number of non-zero entries in a sparse vector, but sampling from the posterior density may not  consistently support the level of sparsity. This was clearly demonstrated by the computed examples with parameter values $r=1$ and $r=1/2$: While the MAP estimates localize well the discontinuities, the draws from the posterior seem to have significantly more discontinuities. This finding underlines also the observation that for sparse recovery, the MAP estimate may be a better summary of the posterior than the posterior mean, which in our example resembles more a smooth reconstruction than a discontinuous one. 

Finally, we point out that the analysis above did not address the question of data sensitivity. It has been demonstrated that if the data have variable sensitivity to different components of the unknown, both MAP and posterior mean estimates may fail to recognize some of the non-zero components in the generative model. While the sensitivity analysis could have been included in the discussion here, it was omitted in order to keep the focus on the sampling techniques proposed in this article, but may be the topic of future work.

\section*{Acknowledgements}

The work of DC was partly supported by the NSF grant DMS 1951446, and that of  ES by the NSF grant DMS-2204618.

\end{document}